\documentclass[11pt]{article}
\usepackage[american]{babel}
\usepackage{times}
\usepackage{amsmath}
\usepackage{latexsym}
\usepackage{amssymb}
\usepackage{theorem}
\usepackage{diagrams}
\theoremstyle{change}
\newtheorem{Thm}{Theorem}[section]
\newtheorem{Cor}[Thm]{Corollary}
\newtheorem{Prop}[Thm]{Proposition}
\newtheorem{Lem}[Thm]{Lemma}
{\theorembodyfont{\rmfamily}
\newtheorem{Num}[Thm]{}

\newtheorem{Def}[Thm]{Definition}}
\newcommand{\proof}{\par\medskip\rm\emph{Proof. }}
\newcommand{\qed}{\ \hglue 0pt plus 1filll $\Box$}

\newcommand{\SKIP}[1]{}

\newcommand{\CC}{\mathbb{C}}   
\newcommand{\HH}{\mathbb{H}}
\newcommand{\HP}{\mathbb{H}\mathrm{P}}

\newcommand{\CP}{\mathbb{C}\mathrm{P}}
\newcommand{\OP}{\mathbb{O}\mathrm{P}}
\newcommand{\RR}{\mathbb{R}}

\newcommand{\QQ}{\mathbb{Q}}

\newcommand{\OO}{\mathbb{O}}

\renewcommand{\SS}{\mathbb{S}}
\newcommand{\ZZ}{\mathbb{Z}}

\newcommand{\SO}{\mathrm{SO}}
\renewcommand{\O}{\mathrm{O}}
\newcommand{\G}{\mathrm{G}}

\newcommand{\BSG}{\mathrm{BSG}}
\newcommand{\BG}{\mathrm{BG}}
\newcommand{\BSO}{\mathrm{BSO}}
\newcommand{\BPL}{\mathrm{BPL}}
\newcommand{\PL}{\mathrm{PL}}

\newcommand{\KOtilde}{\widetilde{KO}}
\newcommand{\KTOPtilde}{\widetilde{KTOP}}
\newcommand{\RTOP}{\mathrm{RTOP}}

\newcommand{\bra}[1]{\langle#1\rangle}
\renewcommand{\tilde}{\widetilde}

\newcommand{\id}{\mathrm{id}}
\newcommand{\fbx}{\mathfrak{x}}
\newcommand{\STOP}{\mathrm{STOP}}
\newcommand{\BSTOP}{\mathrm{BSTOP}}

\newcommand{\BSPL}{\mathrm{BSPL}}
\newcommand{\Spin}{\mathrm{Spin}}
\newcommand{\TOP}{\mathrm{TOP}}
\newcommand{\BTOP}{\mathrm{BTOP}}
\newcommand{\BO}{\mathrm{BO}}
\newcommand{\isom}{\cong}
\newcommand{\homot}{\simeq}
\newcommand{\mapstoo}{\longmapsto}
\newcommand{\too}{\rTo}
\newcommand{\pr}{\mathrm{pr}}
\newcommand{\Ahat}{\,{\widehat {\rlap{$\!\varcal A$}\phantom{W}}\!\!}}
\renewcommand{\L}{\rlap{$\varcal L$}\phantom{W}}

\DeclareMathAlphabet{\varcal}{U}{rsfs}{m}{it}
\begin{document}

\title{\bf Projective planes and their look-alikes}
\author{Linus Kramer\thanks{Supported by a Heisenberg fellowship
by the Deutsche Forschungsgemeinschaft}}

\maketitle
In this paper we classify manifolds which look
like projective planes. More precisely, we consider
$1$-connected closed topological manifolds $M$ with integral
homology
\[
H_\bullet(M)\isom\ZZ^3.
\]
A straight-forward application of Poincar\'e duality shows that
for such a manifold there exists a number $m\geq 2$ such that
$H_k(M)=\ZZ$, for $k=0,m,2m$; in particular, $\dim(M)=2m$ is
even. It follows from Adams' Theorem on the Hopf invariant that
$m$ divides $8$.

We construct a family of topological $2m$-manifolds
$M(\xi)$ which
are Thom spaces of certain topological $\RR^m$-bundles 
(open disk bundles) $\xi$ over the sphere $\SS^m$,
for $m=2,4,8$, and which we call models.
This idea seems to go back to Thom and was exploited further
by Shimada \cite{Shim} and
Eells-Kuiper \cite{EK}. A particular case is worked out in some
detail in Milnor-Stasheff \cite{MS} Ch.~20.
However, these authors used vector bundles
instead of $\RR^m$-bundles. We will see that the non-linearity of
$\RR^m$-bundles yields many more manifolds than the construction by
Eells-Kuiper. In \cite{EK} p.~182, the authors expressed the hope that
\emph{``the given combinatorial examples form a complete set [...]
for $n\neq 4$''}. Our results show that in dimension $n=2m=16$,
their construction missed $27/28$ of the (infinitely many)
combinatorial and topological solutions, while in dimension $n=8$,
they obtained all combinatorial, but only half of the
topological solutions.

Next, we determine certain characteristic classes of our models,
and in particular their rational Pontrjagin classes. Using these
characteristic classes and Wall's surgery sequence, we show that
for $m\neq 2$, every manifold which looks like a projective plane
is homeomorphic to one of our models. Thus, we obtain a complete
homeomorphism classification of all manifolds which are
like projective planes.
Furthermore, we determine which models admit DIFF (or PL)
structures.
The case $m=2$ (so $\dim(M)=4$) is different, but there, we
can apply Freedman's classification of closed $4$-manifolds.

We also classify the homotopy types of $1$-connected
Poincar\'e duality complexes $X$ with $H_\bullet(X)\isom\ZZ^3$. This 
homotopy-theoretic version of our main result (which was already
proved in Eells-Kuiper \cite{EK}) is needed in the
course of the homeomorphism classification; since the methods
here are somewhat different from the rest of the paper, I put it
in an appendix. 

\begin{center}
\textbf{Main results.}
\end{center}
{\em 
Let $M$ be a 1-connected closed topological manifold which looks
like a projective plane, i.e. $H_\bullet(M)\cong\ZZ^3$.
Then $\dim(M)=2m=4,8,16$.

If $m=2$, then $M$ is homeomorphic to the complex projective
plane $\CP^2$ or to the
Chern manifold $\mathrm{Ch}^4$. These two manifolds are
topologically distinguished by their Kirby-Siebenmann numbers
$ks[M]\in\ZZ/2$;
the Chern manifold ($ks[\mathrm{Ch}^4]\neq 0$) admits no DIFF structure
(this case is due to Freedman \cite{Fr}).

If $m=4$, then $M$ is homeomorphic to one of our models
$M(\xi)$. Topologically, it is determined by the
Pontrjagin number $p_4^2[M]\in\{2(1+2t)^2\mid t\in\ZZ\}$
and the Kirby-Siebenmann number
$ks^2[M]\in\ZZ/2$. These data also determine the oriented bordism class
of $M$, so no two models are equivalent under oriented bordism.
The manifold admits
a PL structure (unique up to isotopy) if and only if $ks^2[M]=0$.

If $m=8$, then $M$ is homeomorphic to one of our models $M(\xi)$. 
Topologically, it is determined by the
Pontrjagin number
$p_8^2[M]\in\left\{\frac{36}{49}(1+2t)^2\mid t\in\ZZ\right\}$ 
and a characteristic number
$\frac76p_8\kappa[M]\in\ZZ/4$,
determined by the integral characteristic class
$\frac76p_8(M)$ and a certain $8$-dimensional
PL characteristic class $\kappa$ with $\ZZ/4$-coefficients.
These data also determine the oriented bordism class
of $M$, so no two models are equivalent under oriented bordism.
These manifolds admit a PL structure (unique up to isotopy).}

We determine also which of these manifolds admit a DIFF structure,
and determine the homotopy type in terms of the characteristic classes.
See Sections \ref{ClassifyModels} and \ref{HomotopyOfModels} for
more detailed statements. 
As a by-product of our proof,
we obtain an explicit classification of $\RR^m$-bundles
over $\SS^m$ in terms of characteristic classes, for $m=2,4,8$.

\begin{center}
\textbf{* * *}
\end{center}
Topological geometry plays no r\^ole in this paper. However,
the motivation to write it came from a long-standing
open problem in topological geometry: 
\begin{description}
\item[$(*)$] \emph{What are the possible
homeomorphism types
of the point spaces  of compact projective planes
(in the sense of Salzmann \cite{CPP})?}
\end{description}
The point space $P$ of a compact projective plane
is always the Thom space of a locally compact fiber bundle,
see Salzmann \emph{et al.} \cite{CPP} Ch.5, in particular 51.23
(problem $(*)$ should not be confused with the
\emph{geometric} problem of classifying all compact projective planes
with large automorphism groups which was
solved by Salzmann and his school \cite{CPP}).

Now $(*)$ turns out to be a difficult problem. The present
state of affairs is as follows, see \cite{CPP}.
Let $P$ be the point space
of a compact projective plane.
If the covering dimension
of $P$ is $\dim(P)=0$, then $P$
is either finite or homeomorphic to the Cantor set
$\{0,1\}^{\mathbb{N}}$.
If $1\leq \dim(P)\leq 4$, then
$P\isom\RR\mathrm{P}^2$ or
$P\isom\CC\mathrm{P}^2$;
this was proved by Salzmann and Breitsprecher already
in the late 60s \cite{Bs} (surprisingly, this did not require
results about $4$-manifolds). The proof
depends on a result by Borsuk about low-dimensional ANRs and
on Kneser's Theorem $\SO(2)\homot\STOP(2)$.
In (finite) dimensions bigger than $4$, 
L\"owen \cite{Lo} applied sheaf-theoretic cohomology to the problem.
Using a beautiful local-to-global
argument, he proved that $P$ is an $m-1$-connected
Poincar\'e duality complex and an integral $2m$-dimensional ENR manifold with
$H_\bullet(P)\isom\ZZ^3$, and that $m=2,4,8$.

So the topological problem $(*)$ is reduced to the following steps.
(1) Prove that the topological dimension $\dim(P)$ is finite.
(2) Assuming that $\dim(P)<\infty$, prove that
$P$ is a manifold (and not just an integral ENR manifold).
(3) Assuming that $P$ is a manifold, determine its
homeomorphism type.

Each step seems to be difficult.
Under the additional
assumption that the compact projective plane is \emph{smooth}
(in the sense of \cite{CPP}: the geometric operations are smooth maps),
a complete homeomorphism classification (based on characteristic classes)
of the point spaces was carried out in \cite{KrSmooth}.
Buchanan \cite{Buch} determined the homeomorphism types of the
point spaces of compact projective planes coordinatized by real
division algebras by a direct homotopy-theoretic argument
(note that besides the classical alternative division algebras
$\RR$, $\CC$, $\mathbb H$, and
$\mathbb O$, there exists a continuum of other real division algebras).
In both cases, the homeomorphism types of the point spaces
turn out to be the classical ones, $\RR\mathrm P^2$,
$\CC\mathrm P^2$,
$\mathbb H\mathrm P^2$, or $\mathbb O\mathrm P^2$.
My hope is that the results in this paper, together
with Knarr's Embedding Theorem \cite{Kn} and the result in
\cite{KrP=L}
will eventually lead to a solution of (3). 
\begin{center}
\textbf{* * *}
\end{center}
I have tried to make the paper self-contained and accessible to
non-experts. There is necessarily a certain overlap with the
paper by Eells-Kuiper \cite{EK}. My aim was to give complete
proofs for all steps of the classification, starting only from
general facts about bundles and manifolds. Thus, the reader
is not assumed to be familiar with \cite{EK} (although this
fundamental paper is certainly to be recommended).

\textbf{Standing assumptions.}
An $n$-manifold (without boundary) is a metrizable,
second countable space which is locally homeomorphic to $\RR^n$.
Throughout, all maps are assumed to be continuous.
Except for the appendix, maps and homotopies are not required
to preserve base points, unless the contrary is stated explicitly.
Thus $[X;Y]$ denotes the set of all free homotopy classes of maps
from $X$ to $Y$.
If $X,Y$ are well-pointed spaces, and if $Y$ is $0$-connected,
then the fundamental group $\pi_1(Y)$ acts on the set $[X;Y]_0$
of based homotopy classes; the set $[X;Y]$ of all
free homotopy classes can be identified with the orbit set
of this action, see Whitehead \cite{Whi} Ch.~III.1.
If $Y$ is an $H$-space (or if $Y$ is $1$-connected) this action
is trivial, so $[X;Y]=[X;Y]_0$.

\medskip\noindent
\textbf{Acknowledgements.} I am indebted to Stephan Stolz and
Michael Weiss for sharing some of their insights.

\section{Preliminaries on bundles}

Fiber bundles, microbundles, fibrations, and Thom spaces play a prominent
r\^ole in this paper, so we briefly recall the relevant notions.
We refer to Holm \cite{Ho}, Milnor \cite{MiMicro}, Dold \cite{Dold},
and to the books by Kirby-Siebenmann \cite{KS}, Rudyak \cite{Rud}, and
Husemoller \cite{Hus}.
\begin{Num}
A \emph{bundle} $\phi=(E,B,p)$ over a space $B$ is a
map $E\too^p B$. The class of all bundles over $B$ forms in an
obvious way a category
whose morphisms $\phi\too^f\phi'$ are 
commutative diagrams
\begin{diagram}[size=2em]
E' & \too^f & E\\
\dTo^{p'}&&\dTo_{p}\\
B & \rEq & B.
\end{diagram}
An isomorphism in this category is called an \emph{equivalence} of bundles
and denoted $\phi\cong\phi'$;
in the diagram above, $f$ is an equivalence if and only if $f$ is
a homeomorphism. The categorical product of two bundles $\phi$, $\phi'$
is the
\emph{Whitney sum} $\phi\oplus\phi'$; its total space is
$E\oplus E'=\{(e,e')\in E\times E'\mid p(e)=p'(e')\}$, with the
obvious bundle projection.
A homotopy between two morphisms
$f_0,f_1:\phi'\pile{\too\\\too} \phi$ or \emph{homotopy over $B$}
is a homotopy $E'\times[0,1]\too E$ with the property
that the diagram
\begin{diagram}[size=2em]
E' & \rTo^{f_t} & E \\
\dTo && \dTo \\
B & \rEq & B
\end{diagram}
commutes for all $t\in[0,1]$. A morphism $f$ is called a
\emph{fiber homotopy equivalence} if it has a homotopy inverse
bundle map $g$, i.e. if $fg$ and $gf$ are homotopic over $B$
to the respective identity maps; in this case we write $\phi\simeq\phi'$.
A \emph{section} of a bundle $\phi=(E,B,p)$ is a morphism $s$
from the identity bundle $(B,B,\id_B)$ to $\phi$,
\begin{diagram}[size=2em]
&& E\\
&\ruTo^s&\dTo_{p}\\
B & \rEq & B.
\end{diagram}
and we call $(E,B,p,s)$ a \emph{sectioned bundle}.

A map $B'\too^g B$ induces a contravariant functor $g^*$ which
assigns to every bundle $\phi=(E,B,p)$ the \emph{pull-back bundle}
$g^*\phi=(g^*E,B',p')$, with $g^*E=\{(e,b')\in E\times B'\mid p(e)=g(b')\}$
and $p'(e,b')=b'$. If $g$ is a homeomorphism, then $g^*$ is an equivalence
of categories; in this case, two bundles $\phi$ and $\phi'$ are called
\emph{weakly equivalent} if $\phi'$ is equivalent to $g^*\phi$, in other
words, if there are homeomorphisms
\begin{diagram}[size=2em]
E' & \too^f_\cong & E\\
\dTo^{p'}&&\dTo_{p}\\
B' & \too^g_\cong & B
\end{diagram}
commuting with the bundle projections; such a weak equivalence is
denoted $\phi\cong_g \phi'$.
\end{Num}
\begin{Num}
A bundle is called a \emph{fibration} if the homotopy
extension problem
\begin{diagram}[size=2em]
X\times\{0\}& \rTo^f & E \\
\dInto & \ruDotsto & \dTo_p \\
X\times[0,1] & \rTo^g & B
\end{diagram}
has a solution for every space $X$. We call a fibration
\emph{$n$-spherical} if every fiber $E_b=p^{-1}(b)$ has the
homotopy type of an $n$-sphere.
\end{Num}
For a subspace $A\subseteq B$, we have the \emph{restriction}
$\phi|_A=(E_A=p^{-1}(A),A,p|_{E_A})$ of the bundle $\phi$.
\begin{Def}
A bundle is called
a \emph{fiber bundle} with typical fiber $F$ if every $b\in B$ has
an open neighborhood $U$ such that the restriction
$\phi|_U$ is equivalent to
the product (or trivial) bundle $(F\times U,U,\pr_2)$
\begin{diagram}[size=2em]
F\times U & \too_\cong & E_U & \rInto & E\\
\dTo^{pr_2}&&\dTo && \dTo\\
U & \rEq & U & \rInto & B;
\end{diagram}
such a local trivialization is also called a \emph{coordinate chart}
for the bundle.
If in addition a base point is fixed in the fiber $F$, one obtains
in an obvious way a \emph{sectioned fiber bundle}.
\end{Def}
For technical reasons, it is often convenient to consider numerable
fiber bundles. For example, every numerable fiber bundle is
automatically a fibration, see Spanier \cite{Span} Ch.~2.7 Theorem 12.
\begin{Def}
A locally finite covering $\{V_i\mid i\in I\}$ of $B$ by open sets
is called \emph{numerable} if there
exist maps $f_i:B{\too}[0,1]$ with $f_i^{-1}((0,1])=V_i$,
such that $\sum_{i\in I}f_i=1$. A fiber bundle 
is called \emph{numerable} if there exists a numerable covering of $B$
by coordinate charts. In our setting, most base
spaces will be paracompact,
so fiber bundles are automatically numerable.
\end{Def}
\begin{Def} An \emph{$n$-sphere bundle} is a numerable fiber bundle with
$\SS^n$ as typical fiber. An \emph{$\RR^n$-bundle} is a sectioned numerable
bundle with $(\RR^n,0)$ as typical fiber; the section is denoted
$s_0$ and called the \emph{zero-section}. The trivial $\RR^n$-bundle
(over any space) will be denoted $\underline\RR^n$; its total space is
$E=\RR^n\times B$, with $p=\pr_2$.
An $n$-dimensional \emph{vector bundle} is an $\RR^n$-bundle which
carries in addition a real vector space structure on each fiber which is
compatible with the given coordinate charts.
Two $\RR^n$-bundles (or vector bundles) $\xi$, $\xi'$ are called
\emph{stably equivalent} if there is an equivalence
$\xi\oplus\underline\RR^k\cong\xi'\oplus\underline\RR^{k'}$,
for some $k,k'\geq 0$.
\end{Def}
A crucial property of $\RR^n$-bundles is the following homotopy
property.
\begin{Lem}
\label{HomotopyInvarianceOfBundles}
Let $\xi$ be an $\RR^n$ bundle over $B$.
If $g_0,g_1:B'\pile{\too\\\too}B$ are homotopic, then
there is an equivalence $g_0^*\xi\cong g_1^*\xi$.

\proof See Holm \cite{Ho} Lemma~1.5.
\qed
\end{Lem}
\begin{Num}
\label{ZeroSection}
In an $\RR^n$-bundle $\xi=(E,B,ps_0)$, the
zero-section $s_0$ is a homotopy inverse to the bundle projection
$p$, and $s(B)$ is a strong deformation retract of the total space
$E$, see Holm \cite{Ho} Theorem~3.6. In particular, the section
$s_0:B\too E$ is a cofibration. It follows that the quotient
$E/s_0(B)$ is contractible.
\end{Num}
\begin{Def}
\label{ThomSpaceDef}
From each $\RR^n$-bundle $\xi$, one obtains an $n$-sphere bundle 
$^s\xi$ by compactifying each fiber of $\xi$. The resulting
bundle $^s\xi$ has two sections, the zero-section $s_0$ and the
section $s_\infty$ corresponding to the new points added in the fibers.
Let $E$ denote the total space of $\xi$, and put $E_0=E\setminus s_0(B)$,
the total space with the zero-section removed.
Finally, let $^uE=E_0\cup s_\infty(B)$. Then clearly,
$^u\xi=(^uE,B,^up,s_\infty)$ is again an $\RR^n$-bundle
(called the \emph{upside down bundle} in \cite{KrP=L}),
and $E_0\too B$ is a numerable fiber bundle with $\RR^n\setminus0$
as typical fiber. We call $E_0\too B$ the \emph{spherical fibration}
corresponding to the $\RR^n$-bundle $\xi$.
The \emph{Thom space} $M(\xi)$ of an $\RR^n$-bundle $\xi$ is the quotient
\[
M(\xi)=E\cup s_\infty(B)/s_\infty(B)=E\cup\{o\}.
\]
\end{Def}
We denote the base point (the tip) of $M(\xi)$ by $o$.
If $B$ is compact, then
$M(\xi)=E\cup\{o\}$ is the same as the one-point compactification of $E$.
\begin{Num}
\label{ThomRem}
By the previous remarks, $o$ is a strong deformation retract of
$E_0\cup\{o\}$. In particular, there is a natural (excision) isomorphism
$H^\bullet(E,E_0)\cong H^\bullet(M(\xi),o)$.
\end{Num}
We need also the concept of a microbundle, cp.~Milnor \cite{MiMicro}
and Holm \cite{Ho}.
\begin{Def}
\label{MicroBundleDef}
An \emph{$n$-microbundle} $\fbx=(E,B,p,s)$ is a sectioned bundle,
subject to the following
condition: for every $b\in B$ there exists an open neighborhood
$U$ of $b$, an open subset $V\subseteq E_U=p^{-1}(U)$ containing $s(U)$
and a section-preserving
homeomorphism $h:U\times\RR^n\too V$ such that the diagram
\begin{diagram}[size=2em]
\RR^n\times U & \too^h & V & \rInto & E_U & \rInto & E\\
\dTo^{pr_2}&&\dTo && \dTo && \dTo \\
U & \rEq & U & \rEq & U & \rInto & B;
\end{diagram}
commutes. The difference between a microbundle and an $\RR^n$-bundle is
that $h$ need not be surjective onto $E_U$.
Similarly as for $\RR^n$-bundles, we require the existence
of a numerable covering of $B$ by such local charts.
\end{Def}
Clearly, every $\RR^n$-bundle is an $n$-microbundle.
It is also clear that there exist microbundles which
are not fiber bundles.
A particularly important example is the \emph{tangent microbundle}
$\mathfrak t M$
of a manifold $M$: here, $E=M\times M$, the bundle projection is
$\pr_2$ and the section is the diagonal map, $s(x)=(x,x)$,
see Milnor \cite{MiMicro} Lemma~2.1.

The Kister-Mazur Theorem (see Theorem \ref{KisterMazur} below)
says that
microbundles are in a sense equivalent to $\RR^n$-bundles, a fact
which is not obvious at all.
If $E'\subseteq E$ is a neighborhood of $s(B)$, then it is not
difficult to see that $E'\too B$ is again a microbundle
$\fbx'$ \emph{contained in} $\fbx$. Two
microbundles  $\fbx_1,\fbx_2$ over the same base
$B$ are called \emph{micro-equivalent} if they contain microbundles
$\fbx_1',\fbx_2'$ which are equivalent as bundles
(this is also sometimes called a
micro-isomorphism or an isomorphism germ).
In the case of numerable microbundles one has to be careful: a
microbundle contained in a numerable microbundle need \emph{a priori}
not be numerable.
\begin{Thm}[Kister-Mazur]
\label{KisterMazur}
Let $\fbx$ be a numerable $n$-microbundle. Then there exists
a numerable microbundle $\fbx'$ contained in $\fbx$
which is an $\RR^n$-bundle, and $\fbx'$ is unique up to
equivalence.

\proof See Holm \cite{Ho} Theorem 3.3.
\qed
\end{Thm}
In particular, the tangent microbundle $\mathfrak t M$
of any (metrizable) $n$-manifold
$M$ contains an $\RR^n$-bundle, unique up to equivalence. We will
call this $\RR^n$-bundle $\tau M$ (and any bundle equivalent to it) the
\emph{tangent bundle} of $M$. If $M$ happens to be a smooth manifold,
one can show that $\tau M$ is equivalent to the smooth tangent bundle
$TM$, see Milnor \cite{MiMicro} Theorem~2.2.

If $\xi$ is an $\RR^n$-bundle over a manifold $B$, then the total
space $E$ is clearly a manifold, and the zero-section $s_0(B)$
is a submanifold with normal (micro) bundle $\xi$, see Milnor
\cite{MiMicro} Sec.~5. We require the following splitting result.
\begin{Prop}
\label{SplitsOverSection}
There is an equivalence
\[
s_0^*\tau E\cong\tau B\oplus\xi.
\]

\proof This follows from Milnor \cite{MiMicro} Theorem~5.9, combined
with the Kister-Mazur Theorem~\ref{KisterMazur} above.
\qed
\end{Prop}

\section{Constructing the models}

In this section we construct a family of manifolds as Thom spaces
of $\RR^m$-bundles over the sphere $\SS^m$, for $m=2,4,8$,
which we call \emph{models}.
We begin with some general remarks about Thom spaces of
$\RR^m$-bundles over $\SS^m$. We fix a generator
$[\SS^m]\in H_m(\SS^m)$.
Let $\xi$ be an $\RR^m$-bundle over $\SS^m$, for $m\geq 2$, with
total space $E$, and let
$E_0=E\setminus s_0(\SS^m)$ denote the total space with the zero-section
removed. Since $\SS^m$ is 1-connected, the bundle $\xi$ is
orientable, and we may choose an \emph{orientation class}
$u(\xi)\in H^m(E,E_0)$,
see Spanier \cite{Span} Ch.~5.7 Corollary 20.
The image $e(\xi)=s_0^\bullet(u(\xi)|_E)$ of $u(\xi)$
in $H^m(\SS^m)$ is the \emph{Euler class} of $\xi$.
We call the integer 
\[
|e|=|\bra{e(\xi),[\SS^m]}|
\]
the \emph{absolute Euler number} of $\xi$; it is
independent of the choice of $[\SS^m]$ and of $u(\xi)$.
Let $M(\xi)$ denote the Thom space of $\xi$, and let
\[
\Phi:H^\bullet(B)\too^\cong H^{\bullet+m}(E,E_0),\qquad
\Phi(v)=p^\bullet(v){\smile}u(\xi)
\]
denote the Gysin-Thom isomorphism, see Spanier
\cite{Span} Ch.~5.7 Theorem~10.
By \ref{ThomRem}, this yields -- via excision -- isomorphisms
\[
H^\bullet(\SS^m)\cong H^{\bullet+m}(E,E_0)\cong
H^{\bullet+m}(M(\xi),o).
\]
Let $y_m,y_{2m}$ be generators for the infinite cyclic groups
$H^m(M(\xi))$ and $H^{2m}(M(\xi))$, respectively.
\begin{Lem}
\label{CohomologyOfThomSpace}
In the cohomology ring $H^\bullet(M(\xi))$,
we have the relation $y_m^2=\pm |e| y_{2m}$.

\proof
Let $x\in H^m(\SS^m)$ be the generator dual to $[\SS^m]$
and let $\Phi$ denote the Gysin-Thom isomorphism.
Thus $\Phi(1)=u(\xi)=\pm y_m$ and $\Phi(x)=\pm y_{2m}$. Let
$e(\xi)=\varepsilon x$,
for $\varepsilon\in\ZZ$ (so $|e|=|\varepsilon|$).
Then $\Phi(e(\xi))=\varepsilon\Phi(x)=u(\xi){\smile}u(\xi)$,
since $p^\bullet(e(\xi))=u(\xi)|_E$.
\qed
\end{Lem}
\begin{Num}
If $m$ is odd, then $u(\xi)\smile u(\xi)=0$, so $e(\xi)=0$. Therefore,
$|e|=0$ if $m$ is odd.
\end{Num}
\begin{Prop}
\label{IfManifolde=1}
If $M(\xi)$ is a manifold, then $|e|=1$ and $m$ is even and divides $8$.
Moreover, $H^\bullet(M(\xi);R)\cong R[y_m]/(y_m^3)$ for any commutative
ring $R$.

\proof
If $M(\xi)$ is a manifold with fundamental class $\mu$, then
Poincar\'e duality implies that the map
\[
H^m(M(\xi))\otimes H^m(M(\xi))\too\ZZ,\qquad
u\otimes v\mapstoo\bra{u{\smile}v,\mu}
\]
is a duality pairing, so $|e|=1$ and $m$ is even. Thus
$H^\bullet(M(\xi);R)\cong R[y_m]/(y_m^3)$ for any commutative ring $R$.
Since $M(\xi)$ is a manifold, it is an ANR, see 
Hanner \cite{Han} Theorem~3.3 or Hu \cite{HuRetr} p.~98
and thus homotopy equivalent
to a CW-complex $X$, see Weber \cite{Web} p.~218;
by standard obstruction theory,
$X\simeq\SS^m\cup e^{2m}$
is homotopy equivalent to a $2$-cell complex, see Wall \cite{WallCW}
Proposition~4.1.
By Adams-Atiyah \cite{AA} Theorem~A, this implies that $m=2,4,8$.
More details can be found in the appendix.
\qed
\end{Prop}
The exact homotopy sequence of the $m-1$-spherical fibration
$E_0\too\SS^m$ shows that $\pi_1(E_0)$ is abelian, and that
$E_0$ is $m-2$-connected.
For $|e|=1$, the Gysin sequence
\[
\too H_k(E_0)\too H_k(\SS^m)\too^{e(\xi)\frown\ \ }
H_{k-m}(\SS^m) \too H_{k-1}(E_0)\too
\]
breaks down into isomorphisms
$H_0(E_0)\cong H_0(\SS^m)$ and 
$H_m(\SS^m)\cong H_{2m-1}(E_0)$, and all other homology groups of
$E_0$ vanish. A repeated application of the Hurewicz isomorphism
shows that $E_0$ is $2(m-1)$-connected, with
$\pi_{2m-1}(E_0)\cong H_{2m-1}(E_0)\cong\ZZ$.
Thus
\[
H_k(E_0\cup\{o\},E_0)\cong
\tilde H_{k-1}(E_0)\cong
\begin{cases}\ZZ & \text{ for }k=2m\\
0 & \text{ else,}\end{cases}
\]
because $E_0\cup\{o\}$ is contractible.
In other words, $M(\xi)$ has for $|e|=1$ the same local
homology groups at $o$ as $\RR^{2m}$.

Recall that, due to the embeddability of second countable
finite dimensional metric spaces,
a locally compact finite dimensional second countable
ANR (absolute neighborhood retract for the class of metric spaces
see Hu \cite{HuRetr}) is exactly the
same as an ENR (euclidean neighborhood retract, see 
Hurewicz-Wallman \cite{HuWa} Ch.~V, Engelking \cite{EngDim}
Theorem 1.11.4, 
and Dold \cite{DoldAT} Ch.~IV.8).
\begin{Lem}
For $|e|=1$, the Thom space $M(\xi)=E\cup\{o\}$ is an
integral ENR $2m$-manifold, i.e. an ENR (euclidean neighborhood retract)
which has the same local
homology groups as $\RR^{2m}$.

\proof
The space $E\cup s_\infty(\SS^m)$ is a $2m$-manifold
(and in particular an ENR),
and $M(\xi)=E\cup\{o\}=E\cup s_\infty(\SS^m)/s_\infty(\SS^m)$ is
a quotient of an ENR (the manifold $E\cup s_\infty(\SS^m)$)
by a compact ENR subspace (the $m$-sphere $s_\infty(\SS^m))$.
Such a quotient
is again an ENR, see Hanner \cite{Han} Theorem~8.2 or
Hu \cite{HuRetr} Ch.~IV.
The local homology groups at $o$ were determined
above; every point in $E$ has a locally euclidean neighborhood
and thus the same local homology groups as $\RR^{2m}$.
\qed
\end{Lem}
Our next aim is to show that $M(\xi)$ is in fact a manifold. Since
$M(\xi)=E\cup\{o\}$ and $E$ is a manifold,
the only point which we have to
consider in detail is $o$. First, we prove that $o$ is
$2(m-1)$-LC in $M(\xi)$, i.e. that every open neighborhood
$V$ of $o$ contains an open neighborhood $V'$ of $o$
such that for $k\leq 2(m-1)$, every map $\SS^k\too V'\setminus\{o\}$
is homotopic in $V\setminus\{o\}$ to a constant map.
Clearly, we are done if we can show that $V'\setminus\{o\}$
is $2(m-1)$-connected.
\begin{Lem}
The space $M(\xi)$ is $2(m-1)$-LC at $o$ if $|e|=1$.

\proof
Let $^uE=E_0\cup s_\infty(\SS^m)$ denote the upside-down bundle
obtained from $\xi$, cp.~Definition~\ref{ThomSpaceDef}.
Then $^uE\too\SS^m$ is an
$\RR^m$ bundle and in particular a microbundle,
cp.~Definition~\ref{MicroBundleDef}.
Let $V$ be an open neighborhood of $o$ in
$M(\xi)$, and let $f$ denote the map
$^uE\too E_0\cup\{o\}\subseteq M(\xi)$ which
collapses the $s_\infty$-section to
the point $o$.
Then $U=f^{-1}(V)$ is an open neighborhood of $s_\infty(\SS^m)$
in the upside-down bundle $^uE$.
By the Kister-Mazur Theorem \ref{KisterMazur},
there exists an open neighborhood $U'$ of $s_\infty(\SS^m)$ in
$U$ with $^uE\supseteq U \supseteq U'$, such that
$U'\too\SS^m$ is an $\RR^m$-bundle
equivalent to $^uE\too\SS^m$. In particular,
$U'\setminus s_\infty(\SS^m)\cong E_0$ is $2(m-1)$-connected.
Now we put $V'=f(U')$.
\qed
\end{Lem}

\begin{Cor}
For $|e|=1$, the Thom space $M(\xi)$ is a $1$-connected closed
$2m$-manifold.

\proof
The space $M(\xi)\setminus\{o\}$ is a $2m$-manifold, and
$M(\xi)$ is $1$-LC at $o$. Thus, $o$ has an open neighborhood
homeomorphic to $\RR^{2m}$; for $m=2$, this follows from
Freedman-Quinn \cite{FQ} Theorem~9.3A
(and also from Kneser's Theorem $\TOP(2)\simeq\O(2)$,
see Theorem \ref{TheBundlesWeObtain} below),
and for $m=4,8$ from Quinn \cite{Qui} Theorem~3.4.1.
Van Kampen's Theorem, applied to the diagram
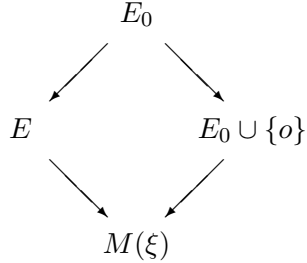
\begin{diagram}[size=2em]
& & E_0 \\
 & \ldTo & & \rdTo \\
E &&& &E_0\cup\{o\} \\
 & \rdTo && \ldTo \\
   && M(\xi)
\end{diagram}
shows that $\pi_1(M(\xi))=0$,
because $E$, $E_0\cup\{o\}$, and $E_0$ are $1$-connected.
\qed
\end{Cor}
\begin{Def}
The manifolds $M(\xi)$ obtained in this way as Thom spaces of
$\RR^m$-bundles with $|e|=1$ will be called \emph{models}.
\end{Def}
The same argument as above
shows that $E/s_0(\SS^m)$ is a manifold, and so
$S=M(\xi)/s_0(\SS^m)$ is a manifold, too. Similarly as above,
Van Kampen's Theorem shows that $S$ is $1$-connected.
As $E_0$ has the same homology as $\SS^{2m-1}$,
the exact homology sequence of the pair $(E,E_0)$ shows that
the composite $\SS^m\too^{s_0}E\too(E,E_0)$ is an isomorphism
in $m$-dimensional homology. Thus $\tilde H_\bullet(S)\cong
H_\bullet(M(\xi),s_0(\SS^m))$ (here we use that
$s_0:\SS^m\too M(\xi)$ is a cofibration). Therefore, $S$ is
a 1-connected homology $2m$-sphere, and thus, by the proof
of the generalized Poincar\'e conjecture in higher dimensions,
a sphere (in dimension $4$, see Freedman \cite{Fr}, and in higher
dimension Smale \cite{Sma} and Newman \cite{New}).
In particular, $E_0\cup\{o\}\cong\RR^{2m}$ is an open cell.
Thus, $M(\xi)$ is a compactification of an open $2m$-cell
by an $m$-sphere.
\begin{Prop}
Each model $M(\xi)$ can be decomposed as 
\[
M(\xi)=X\mathrel{\dot\cup} U,
\]
with $X=s_0(\SS^m)$ homeomorphic to $\SS^m$, and $U=E_0\cup\{o\}$
open, dense,
and homeomorphic to $\RR^{2m}$.
\qed
\end{Prop}

\section{Homeomorphisms between different models}
In the last section, we constructed for every
oriented $\RR^m$-bundle $\xi$ over $\SS^m$ with absolute Euler number
$|e|=|\bra{e(\xi),[\SS^m]}|=1$ a manifold $M(\xi)$. In this
section, we determine under which conditions
there are homeomorphisms $M(\xi)\cong M(\xi')$ between different
models.
Clearly, a weak bundle equivalence
$\xi\cong_g \xi'$ induces a homeomorphism $M(\xi)\cong M(\xi')$
between the Thom spaces. We will see that this is in fact
the only possibility. In the first part of this section,
we assume only that $m$ is even; in the second part, we return to
the special case of our models where $m=2,4,8$ and $|e|=1$.

Let $\xi$ be an $\RR^m$-bundle
over $\SS^m$, with absolute Euler number $|e|$, for $m\geq 2$
even. Let $X=s_0(\SS^m)\subseteq E$.
By Proposition \ref{SplitsOverSection}, the tangent bundle $\tau E$
of the manifold $E$ splits along $X$ as a sum of a horizontal
and a vertical bundle.
Since $\tau\SS^m\oplus\underline\RR\cong\underline\RR^{m+1}$, we have
the following result.
\begin{Lem}
Let $\tau E$ denote the topological tangent bundle of
$E$, and let $X=s_0(\SS^m)$. Then $X\subseteq E$ is an embedded
submanifold with normal bundle weakly equivalent to $\xi$ and
\[
\xi\oplus\tau\SS^m\cong s_0^*\tau E.
\]
In particular, $\xi$ and $\tau E|_X$ are weakly stably equivalent,
\[
\xi\oplus\underline\RR^{m+1}\cong s_0^*\tau E\oplus\underline\RR.
\]
\qed
\end{Lem}
Suppose now that there is a homeomorphism
$E\too^f_\cong E'$ of total spaces of $\RR^m$ bundles $\xi,\xi'$
over $\SS^m$, for $m\geq 2$ even. Both 
$s_0$ and $f^{-1}s_0'$ represent generators of $\pi_m(E)\cong\ZZ$.
Thus, there exists a homeomorphism $g:\SS^m\too\SS^m$ of degree
$\pm1$ such that the diagram
\begin{diagram}
E & \rTo^f_\cong & E' \\
\uTo^{s_0} && \uTo_{s_0'} \\
\SS^m & \rTo^g_\cong &\SS^m
\end{diagram}
is homotopy commutative. By Lemma \ref{HomotopyInvarianceOfBundles},
there is an equivalence $s_0^*\tau E\cong g^*(s_0')^*\tau E'$,
whence
$\xi\oplus\underline\RR^{m+1}\cong g^*\xi'\oplus\underline\RR^{m+1}$.
In other words, the bundles $\xi$ and $\xi'$ are weakly stably equivalent.
\begin{Lem}
\label{ThomSpaceHomeo}
Suppose that there is a homeomorphism of Thom spaces $M(\xi)\cong M(\xi')$.
Then there is a homeomorphism between the total spaces $E\cong E'$,
and the bundles $\xi$, $\xi'$ have the same absolute Euler number.

\proof
By Lemma~\ref{CohomologyOfThomSpace},
the absolute Euler number can be seen from the cohomology ring
of $M(\xi)$, so $|e|=|e'|$. If $M(\xi)$ and $M(\xi')$ are manifolds,
then they are homogeneous and the existence of a homeomorphism
$M(\xi)\cong M(\xi')$ implies the existence of a homeomorphism
$E\cong E'$. If $M(\xi)$ and $M(\xi')$ are not manifolds, then
a homeomorphism maps the unique non-manifold point $o$ of $M(\xi)$
onto the unique non-manifold point $o'$ of $M(\xi')$,
and so it maps $E$ onto $E'$.
\qed
\end{Lem}
The proof of the next proposition involves classifying spaces, so we
postpone it to \ref{StabilizingPropProof}.
\begin{Prop}
\label{StabilizingProp}
Let $\xi,\xi'$ be $\RR^m$-bundles over $\SS^m$, for $m\geq 2$ even.
Suppose that there is a stable equivalence
$\xi\oplus\underline\RR^k\cong\xi'\oplus\underline\RR^k$,
and that the absolute Euler numbers of $\xi$ and $\xi'$
are equal, $|e|=|e'|$. Then there is an equivalence $\xi\cong\xi'$.
\end{Prop}
Combining these results, we obtain a complete homeomorphism classification
of the Thom spaces $M(\xi)$, for $m\geq 2$ even, in terms of bundles.
\begin{Prop}
\label{HomeoProp}
Let $m\geq 2$ be even, let $\xi,\xi'$ be $\RR^m$-bundles over $\SS^m$.
If there is a homeomorphism between the Thom spaces
$M(\xi)\cong M(\xi')$, then $\xi$ and $\xi'$ are weakly equivalent.

\proof
By Lemma~\ref{ThomSpaceHomeo} above, the total spaces
$E,E'$ are homeomorphic, and
$|e|=|e'|$. The remarks at the begin of this section show that there
is a weak equivalence between $\xi\oplus\underline\RR^{m+1}$ and
$\xi'\oplus\underline\RR^{m+1}$, induced by a homeomorphism
$g:\SS^m\too\SS^m$. If $\deg(g)=1$, then $g$ is homotopic to the
identity, whence
$\xi\oplus\underline\RR^{m+1}\cong\xi'\oplus\underline\RR^{m+1}$.
By Proposition \ref{StabilizingProp},
this implies that there is an equivalence
$\xi\cong\xi'$.
Otherwise, $\deg(g)=-1$ and we put $\xi''= g^*\xi'$.
Then we have an equivalence
$\xi\oplus\underline\RR^{m+1}\cong\xi''\oplus\underline\RR^{m+1}$,
and, again by Proposition \ref{StabilizingProp}, an equivalence
$\xi\cong\xi''$. But $\xi'$ is weakly equivalent to $\xi''$.
\qed
\end{Prop}
Thus we have reduced the homeomorphism classification of our
models to a classification of $\RR^m$-bundles over $\SS^m$. For the
specific values $m=2,4,8$, this classification will be carried out in
the next section.
We end this section with some simple remarks about characteristic classes
of our models.
Each model $M(\xi)$ has a distinguished orientation: if
$y_m$ is any generator of $H^m(M(\xi))$, then $y_m^2$ is a generator
of $H^{2m}(M(\xi))$ which does not depend on the choice of $y_m$.
So we choose for our models
the fundamental class $[M(\xi)]\in H^{2m}(M(\xi))$
in such a way that
\[
\bra{y_m^2,[M(\xi)]}=1.
\]
Obviously, the Euler characteristic of any model $M(\xi)$ is
\[
\chi_{M(\xi)}=3.
\]
Note also that any homeomorphism $c:\SS^m\too\SS^m$ of degree $-1$
induces a homeomorphism $f_c$ of $M(\xi)$ with the property that
$f_c^\bullet y_m=-y_m$. Thus every graded automorphism of
the cohomology ring $H^\bullet(M(\xi))\cong\ZZ[y_m]/(y_m^3)$
is induced by a homeomorphism.

Recall that the \emph{Wu classes} $v_i\in H^i(M;\ZZ/2)$ of a closed manifold
$M$ are defined by $\bra{v_i\smile x,[M]}=\bra{\mathrm{Sq}^ix,[M]}$.
\begin{Lem}
The total Stiefel-Whitney class of any model $M(\xi)$ is given by
\[
w(M(\xi))=1+y_m+y_m^2,
\]
where $y_m\in H^m(M(\xi);\ZZ/2)$ is a generator. Thus the minimal
codimension for an embedding of $M(\xi)$ in $\SS^{2m+k}$ or $\RR^{2m+k}$
is $k=m+1$.

\proof
We have $\mathrm{Sq}^m\,y_m=y_m^2$, so the total Wu class of $M(\xi)$ is
$v=1+y_m+y_m^2$, and the total Stiefel-Whitney class is
$w(M(\xi))=\mathrm{Sq}\,v=1+y_m+y_m^2$, cp.~Spanier \cite{Span} Ch.~6.10
Theorem~7
and 6.10 8. The non-embedding result follows as in Spanier \cite{Span}
Ch.~6.10 24.
\qed
\end{Lem}
Recall from \ref{ZeroSection} that
the composite $\SS^m\too^{s_0} X\subseteq E$ is a
homotopy equivalence. Since $E_0$ is $2(m-1)$-connected,
the exact homology sequence of the pair $(E,E_0)$ shows that
$H_m(E)\rTo H_m(E,E_0)$ is an isomorphism.
\begin{Lem}
The map $s_0:\SS^m\too M(\xi)$ induces an isomorphism
on the homotopy and (co)ho\-mo\-logy groups up to (and including)
dimension $m$.

\proof
Since $M(\xi)$ is $1$-connected, the claim on the homotopy groups
follows from the corresponding result for the homology groups
by the Hurewicz isomorphism.
\qed
\end{Lem}
\begin{Num}
\label{StableClassInHalfDim}
In particular, if $\gamma(M)$ is a stable $m$-dimensional characteristic
class of the tangent bundle of $M(\xi)$ then
$s_0^\bullet(\gamma(M))=\gamma(\xi)$. For $m=4,8$, this applies
in particular to
the $m$-dimensional rational Pontrjagin class $p_m$ of (the tangent bundle
of) $M(\xi)$,
\[
p_m(\xi)=s_0^\bullet p_m(M(\xi)),
\]
so the $m$-dimensional
Pontrjagin class of $\xi$ determines the $m$-dimensional Pontrjagin class of
$M(\xi)$. A similar result holds for the exotic classes,
the Kirby-Siebenmann class $ks$ and the class $\kappa$ of $M$,
which are constructed in \ref{ExoticClasses}.
This will be used in the Section \ref{CharClass}.
\end{Num}

\section{Characteristic classes}
\label{CharClass}
To get further, we need some results about classifying spaces.
We refer to the books by Milnor-Stasheff \cite{MS},
Kirby-Siebenmann \cite{KS}, Madsen-Milgram \cite{MM} and
to Ch.~IV in Rudyak \cite{Rud}.
We denote the orthogonal group by $\O(n)$, and by $\TOP(n)$ 
the group of all base-point preserving homeomorphisms of $\RR^n$.
The corresponding classifying spaces are $\BO(n)$ and $\BTOP(n)$.
If $X$ is any space, 
then the set of free homotopy classes $[X;\BO(n)]$ is in one to one
correspondence with the equivalence classes of numerable $n$-dimensional
vector bundles over $X$. Similarly, $\BTOP(n)$ classifies 
numerable $\RR^n$-bundles over $X$.
Taking the limit $n\gg1$ one obtains
stable versions $\O$, $\TOP$ of these groups; there are
corresponding classifying spaces $\BO$ and $\BTOP$ which
classify $\RR^n$-bundles up to stable equivalence. 
We need two more classifying spaces. The space $\BPL(n)$ classifies
$\RR^n$-bundles (over simplicial complexes) which admit
PL (piecewise linear) coordinate charts. To us,
the main purpose of $\BPL(n)$
and its stable version $\BPL$ will be the fact that it lies somewhat
in the middle between $\BTOP$ and $\BO$.
Finally, let $\G(n)$ denote the semigroup of all self-equivalences of
the sphere $\SS^{n-1}$. The corresponding classifying space $\BG(n)$
classifies $n-1$-spherical fibrations up to fiber homotopy equivalence.
There are $1$-connected coverings of these spaces
which classify oriented bundles and fibrations: for example, 
the classifying
space $\BSO(n)$ classifies $n$-dimensional oriented vector bundles.
For $k\geq 0$ there are ladders of maps 
\begin{diagram}[size=1.5em]
&& \BO(n) && \too && \BO(n+k) && \too && \BO \\
& \ruTo & \vLine &&& \ruTo &\vLine &&& \ruTo & \dTo \\
\BSO(n) & \HonV & \too && \BSO(n+k) & \HonV & \too && \BSO \\ 
\dTo && \dTo && \dTo && \dTo && \dTo \\
&& \BPL(n) & \hLine & \VonH & \too & \BPL(n+k) &
\hLine & \VonH & \too & \BPL \\
& \ruTo & \vLine &&& \ruTo & \vLine &&& \ruTo \\
\BSPL(n) & \HonV & \too && \BSPL(n+k) & \HonV & \too && \BSPL && \dTo \\
\dTo && \dTo && \dTo && \dTo && \dTo \\
&& \BTOP(n) & \hLine & \VonH & \too & \BTOP(n+k) & \hLine & \VonH &
\too& \BTOP  \\
& \ruTo & \vLine &&& \ruTo &\vLine &&& \ruTo & \dTo \\
\BSTOP(n) & \HonV & \too && \BSTOP(n+k) & \HonV & \too && \BSTOP \\ 
\dTo && \dTo && \dTo && \dTo && \dTo \\
&& \BG(n) & \hLine & \VonH & \too & \BG(n+k) &
\hLine & \VonH & \too & \BG \\
& \ruTo &&&& \ruTo &&&& \ruTo \\
\BSG(n) && \too && \BSG(n+k) && \too && \BSG
\end{diagram}
such that the diagram commutes (at least up to homotopy).
The horizontal arrows correspond to the process of stabilization,
i.e. if $f:X\too\BSTOP(n)$ classifies $\xi$, then the composite
\[
X\too\BSTOP(n)\too\BSTOP(n+k)
\]
classifies $\xi\oplus\underline\RR^k$.
The vertical and the slanted arrows are 'forgetful': they forget the
vector space
structure, the PL structure, and the fiber bundle structure,
respectively, and the slanted arrows forget the orientation.
For these results, see eg.~Rudyak \cite{Rud} Ch.~IV.

By a well-known construction, every continuous map $f:X\too Y$
between topological spaces can be converted into a fibration
$f':P_f\rTo Y$, with $X\simeq P_f$,
see Spanier \cite{Span} Ch.~2.8 Theorem 9.
If $Y$ is path-connected, then all fibers of $f'$ have
the same homotopy type, and it makes sense to speak about the homotopy
fiber of the map $f$.
Let $\mathrm BH\too\mathrm BG$ be one of the maps in the diagram
above. The homotopy fiber of this map is denoted by $G/H$.
One can show that the homotopy fiber of the map $\BO(n)\too\BO(n+k)$
is homotopy equivalent to the Stiefel manifold $\O(n+k)/\O(n)$,
so this terminology fits together with the standard Lie group terminology.

From the homotopy viewpoint, the stable classifying spaces are much
easier to understand. This is partly due to the fact that they
are $H$-spaces (the Whitney sum of bundles is the multiplication).
The homotopy groups of $\BO$ are known by Bott periodicity,
$\BO\times\ZZ\simeq \Omega^8\BO$, see Bott \cite{Bott1} \cite{Bott2}, i.e.
\[
\pi_k(\BO\times\ZZ)\cong\begin{cases}
\ZZ &\text{ for }k\equiv 0\phantom{,2}\pmod 4 \\
\ZZ/2& \text{ for }k\equiv 1,2\pmod 8\\
0 &\text{ else.}
\end{cases}
\]
The homotopy groups of $\BG$ correspond to the stable homotopy groups
of spheres, see eg.~Milnor \cite{MilSpher} \S 2 and
Madsen-Milgram \cite{MM} Ch.~3,
\[
\pi_{k+1}(\BG)\cong\lim_{n\to\infty}\pi_{k+n}(\SS^n)=\pi_{k}^s(\SS^0)
\]
which are known in low dimensions, see eg.~Toda \cite{Toda} Ch.~XIV,
Hu \cite{Hu} pp.~328--332, or
Fomenko-Fuchs-Gutenmacher
\cite{FFG} pp.~300--301.
For the other homotopy fibers, we use the following results
which are obtained from surgery theory.
\begin{Thm}
\label{SurgeryGroups}
The homotopy groups of $\G/\TOP$ are given by the periodicity
$\G/\TOP\times\ZZ\simeq\Omega^4(\G/\TOP)$, see Kirby-Siebenmann
\cite{KS} p.~327 (the $\ZZ$-factor is forgotten there), and
\[
\pi_k(\G/\TOP\times\ZZ)\cong\begin{cases}
\ZZ &\text{ for }k\equiv 0\pmod 4 \\
\ZZ/2& \text{ for }k\equiv 2\pmod 4\\
0 &\text{ else,}
\end{cases}
\]
see also Madsen-Milgram \cite{MM} Ch.~2.
\qed
\end{Thm}
Finally, we use the following result, see Kirby-Siebenmann \cite{KS}
p.~200.
\begin{Thm}
\label{Thetafinite}
The homotopy groups of $\TOP/\O$ and $\TOP/\PL$ are finite in all
dimensions. If $i\geq 5$, then $\pi_i(\TOP/\O)$ is isomorphic to the
Kervaire-Milnor group $\Theta_i$ of DIFF structures on $\SS^i$,
cp.~Kirby-Siebenmann \cite{KS} p.~200, 251, and
$\TOP/\PL$ is an Eilenberg-MacLane space of type $K(\ZZ/2,3)$.
\qed
\end{Thm}
It follows from Serre's $\mathcal C$-theory that
$\tilde H^\bullet(\TOP/\O;\QQ)\cong\tilde H^\bullet(\TOP/\PL;\QQ)
\cong \tilde H^\bullet(\PL/\O;\QQ)\cong 0$, cp.~Spanier \cite{Span}
Ch.~9.6. The Serre spectral
sequence yields thus natural isomorphisms
\[
H^\bullet(\BO;\QQ)\lTo^\cong
H^\bullet(\BPL;\QQ)\lTo^\cong
H^\bullet(\BTOP;\QQ),
\]
see also Kahn \cite{Kahn}.
The rational cohomology ring of $\BO$ is known to be
a polynomial algebra generated by the
\emph{universal Pontrjagin classes} $p_4,p_8,\ldots$,
\[
H^\bullet(\BO;\QQ)\cong\QQ[p_4,p_8,p_{12},\ldots],
\]
see eg.~Madsen-Milgram \cite{MM} p.~13.
In view of the isomorphisms above, there are natural rational
Pontrjagin classes defined for $\RR^n$-bundles and for $\PL$-bundles:
if $X\too\BTOP(n)$ is a classifying map for an $\RR^n$-bundle
$\xi$, then $p_{4k}(\xi)\in H^{4k}(X;\QQ)$
is by definition the pull-back of
the universal Pontrjagin class $p_{4k}\in H^{4k}(\BTOP;\QQ)$
via the composite $X\too\BTOP(n)\too\BTOP$, see also
Milnor-Stasheff \cite{MS} pp.~250--251.

\begin{Num}
Recall that the \emph{signature} of an oriented 
$4k$-manifold $M$ is
by definition the signature of the quadratic form 
$H^{2k}(M;\QQ)\too\QQ$, $v\mapstoo \bra{v^2,[M]}$.
(The signature of a rational quadratic form
represented by a symmetric matrix is the number of strictly positive
eigenvalues minus the number of strictly negative real eigenvalues
of the matrix.)
The signature 
is invariant under oriented bordism and induces a group
homomorphism from the oriented DIFF bordism ring into the
integers,
\[
\mathrm{Sig}:\Omega_\bullet^{\mathrm{SO}}\too\ZZ.
\]
It follows that the signature of a smooth closed $4k$-manifold
can be expressed as a certain
linear combination of the rational Pontrjagin numbers of $M$ which is
given by Hirzebruch's $\L$-genus,
\[
\mathrm{Sig}(M)=\bra{\L_{4k}(M),[M]},
\]
see Milnor-Stasheff \cite{MS} p.~224 and Madsen-Milgram
\cite{MM} Theorem 1.38.
The \emph{Hirzebruch classes} $\L_{4k}\in H^{4k}(\BO;\QQ)$ are certain
rational homogeneous polynomials in the rational Pontrjagin classes,
see Hirzebruch \cite{HirNeue} 1.5, Milnor-Stasheff \cite{MS} \S 19.
These polynomials can be obtained by a formal process from the
power series expansion of the function $f(t)=\frac{\sqrt t}{\tanh\sqrt t}$.
These results were proved by Hirzebruch for smooth closed oriented
manifolds, see Hirzebruch \cite{HirNeue} Hauptsatz 8.2.2.
However, there is an isomorphism
\[
\Omega_\bullet^{\mathrm{SO}}\otimes\QQ\cong
\Omega_\bullet^{\mathrm{STOP}}\otimes\QQ,
\]
so the Signature Theorem carries over to the bordism ring
of closed oriented manifolds, see Kahn \cite{Kahn},
Kirby-Siebenmann \cite{KS} p.~322.
\end{Num}
\begin{Thm}
Let $M$ be a closed oriented $4k$-manifold with fundamental class $[M]$.
Then
\[
\mathrm{Sig}(M)=\bra{\L_{4k}(M),[M]}.
\]
\qed
\end{Thm}
The signature of our models is clearly $\mathrm{Sig}(M(\xi))=1$.
\begin{Lem}
For our models, $M(\xi)$, we have the following relation for the Hirzebruch
classes:
\[
\L_{2m}(M(\xi))=y_m^2.
\]
\qed
\end{Lem}
To compute the Hirzebruch classes in terms of the Pontrjagin classes,
we use the following result which is a direct consequence of
Hirzebruch \cite{HirNeue} 1.4.
Let $f(t)=1+\sum_{k\geq 1}f_kt^k$ be a formal power series, and
let $\{K_k(\sigma_1,\ldots,\sigma_k)\}_{k=1}^\infty$ denote the corresponding
multiplicative sequence, see Hirzebruch \cite{Hir} \S 1.
The associated genus $\varcal F$ of $f(t)$ is defined by
\[
\varcal F_{4k}=K_k(p_4,\ldots,p_{4k})\in H^{4k}(\BO;\QQ).
\]
Hirzebruch's $\L$-genus comes from the formal power series
\[
\ell(t)=\frac{\sqrt t}{\tanh\sqrt t}=
1+\frac13t-\frac1{45}t^2+\frac2{945}t^3-\frac1{4725}t^4+\cdots,
\]
and the $\Ahat$-genus (which will be needed later in
Section \ref{ClassifyModels}) from 
\[
\hat a(t)=\frac{\sqrt t/2}{\sinh\sqrt t/2}=
1-\frac1{24}t+\frac7{5760}t^2-\frac{31}{967680}t^3+
\frac{127}{154828800}t^4+\cdots.
\]
Suppose now that $\xi$ is an $\RR^n$ bundle over a space $X$,
and that $p_{4k}(\xi)$ and $p_{8k}(\xi)$ are the only
non-zero Pontrjagin classes of $\xi$. This holds for example if
$X$ is a space with $H^j(X;\QQ)=0$ for all
$j\neq 0,4k,8k$ (such as our models $M(\xi)$). In the cohomology
ring of such a space, there is the following general formula for 
$\varcal F_{8k}=K_{2k}(0,\cdots,0,p_{4k},0,\dots,0,p_{8k})$.
Put 
\[
f^\vee(t)=f(t)\frac{d}{dt}\left(\frac t{f(t)}\right)=\sum(-1)^ks_kt^k.
\]
By \emph{loc.cit.} 1.4, one has the general relation
\[
K_{2k}(0,\ldots,0,p_{4k},0,\ldots,p_{8k})=s_{2k}p_{8k}
+\frac12(s_k^2-s_{2k})p_{4k}^2
\]
For $\ell^\vee$ and $\hat a^\vee$ one obtains
\begin{align*}
\ell^\vee(t)&=
1-\frac13t+\frac7{45}t^2-\frac{62}{945}t^3+\frac{127}{4725}t^4+\cdots\\
\hat a^\vee(t)&=
1+\frac1{24}t-\frac1{1440}t^2+\frac1{60480}t^3-
\frac1{2419200}t^4+\cdots.
\end{align*}
\begin{Num}
\label{SpecificValues}
For specific values, calculations can readily be done with 
the formal power series package of {\sc Maple}.
In low dimensions, one obtains the following.
\begin{align*}
\L_4(p_4)&=\frac13p_4\\
\Ahat_4(p_4)&=\frac{-1}{2^3\cdot 3}p_4\\
\L_8(p_4,p_8)&=\frac1{3^2\cdot 5}(7p_8-p_4^2)\\
\Ahat_8(p_4,p_8)&=\frac{-1}{2^7\cdot 3^2\cdot 5}(4p_8-7p_4^2)\\
\L_{16}(0,p_8,0,p_{16})&=\frac1{3^4\cdot5^2\cdot 7}(381p_{16}
-19p_8^2)\\
\Ahat_{16}(0,p_8,0,p_{16})&=\frac{-1}{2^{11}\cdot 3^4\cdot5^2\cdot 7}
(12p_{16}-13p_8^2).
\end{align*}
In the last two equations, we assume thus that $p_4=0=p_{12}$.
See also Hirzebruch \cite{HirNeue} 1.5 and 1.6,
Milnor-Stasheff \cite{MS} p.~225,
Lawson-Michelsohn \cite{LM} pp.~231--232,
Eells-Kuiper \cite{EK2} p.~105 for explicit formulas for these
genera in low dimensions.
\end{Num}
From the signature theorem, we obtain thus for our models strong
relations between $p_m$ and $p_{2m}$.
\begin{Lem}
\label{RelationsForModelsLem}
For our models, we have
\begin{align*}
p_4(M)&=\frac13 y_2^2 &&(m=2)\\
p_8(M)&=\frac17(45 y_4^2+p_4(M)^2)&&(m=4)\\
p_{16}(M)&=\frac1{381}(14175 y_8^2+19p_8(M)^2)&&(m=8)
\end{align*}
\qed
\end{Lem}
\begin{Num}
\label{exotic}
We need also the first exotic characteristic classes for
TOP- and PL-bundles. By a standard mapping cylinder construction,
we may convert the maps $\BO\too\BPL\too\BTOP$ into cofibrations,
see Spanier \cite{Span} Ch.~3.2 Theorem 12.
Thus, it makes sense to speak of the homology and homotopy groups
of the pairs $(\BTOP,\BO)$ etc. Note also that these three spaces
are $H$-spaces with isomorphic fundamental groups.
Thus, $\pi_1(\BO)\cong\pi_1(\BPL)\cong\pi_1(\BTOP)\cong\ZZ/2$ acts
trivially on the homotopy groups of each of these pairs,
whence $\pi_k=\pi_k'$ for these pairs
(recall that $\pi_k'$ is the $k$th homotopy group, factored by the action
of $\pi_1$, see Spanier \cite{Span} p.~390). Consequently,
we have Hurewicz isomorphisms $H_k(\BPL,\BO)\cong\pi_k(\BPL,\BO)$
up to and including the lowest dimensions where the right-hand side
is nontrivial, see Spanier \cite{Span} Ch.~7.5 Theorem 4.
Now there is an isomorphism
$\pi_{k-1}(\PL/\O)\cong\pi_k(\BPL,\BO)$, see Whitehead \cite{Whi}
Ch.~IV 8.20. In fact, there is a commutative diagram
\begin{diagram}
\pi_{k-1}(\PL/\O) & \lTo^\partial_\cong & \pi_k(C_{\PL/\O},\PL/\O) &
\rTo_\cong & \pi_k(\BPL,\BO)\\
\dTo && \dTo && \dTo\\
H_{k-1}(\PL/\O) & \lTo^\partial_\cong & H_k(C_{\PL/\O},\PL/\O) &
\rTo & H_k(\BPL,\BO)
\end{diagram}
where $C_{\PL/\O}\simeq *$ is the unreduced cone.
Since $\PL/\O$ is $6$-connected, see Madsen-Milgram p.~33, this diagram
consists of isomorphisms for $k\leq 8$. We obtain a similar diagram
for $\TOP/\O$ and $\TOP/\PL$, with isomorphisms for $k\leq 4$, using
the fact that $\TOP/\PL$ and $\TOP/\O$ are $2$-connected, see
Kirby-Siebenmann p.~246.
\end{Num}
We combine this with the following result due to Hirsch \cite{Hir} p.~356.
\begin{Prop}
There are short exact sequences
\[
0\rTo\pi_k(\BO)\rTo\pi_k(\BPL)\rTo\pi_k(\BPL,\BO)\rTo 0
\]
for all $k\geq 0$. 
\qed
\end{Prop}
\begin{Num}
\label{HirschExact}
Since $\TOP/\PL$ is an Eilenberg-MacLane space of type $K(\ZZ/2,3)$,
and since $\PL/\O$ is $6$-connected, this implies readily that each arrow
\[
\pi_k(\BO)\rTo\pi_k(\BPL)\rTo\pi_k(\BTOP)
\]
is an injection, for all $k\geq 0$.
\end{Num}
We combine this with the
commutative diagram in \ref{exotic} and obtain thus the
following result.
\begin{Lem}
There are natural isomorphisms
$H_k(\BO)\rTo^\cong H_k(\BPL)\rTo^\cong H_k(\BTOP)$ for $k\leq 3$,
and $H_k(\BO)\rTo^\cong H_k(\BPL)$ for $k\leq 7$, and exact
sequences
\begin{diagram}
H_4(\BO) & \rTo & H_4(\BTOP) & \rTo & H_4(\BTOP,\BO) & \rTo & 0\\
H_4(\BPL) & \rTo & H_4(\BTOP) & \rTo & H_4(\BTOP,\BPL) & \rTo & 0\\
H_8(\BO) & \rTo & H_8(\BPL) & \rTo & H_8(\BPL,\BO) & \rTo & 0.
\end{diagram}
Furthermore, $H_4(\BTOP,\BO)\cong H_4(\BTOP,\BPL)\cong\ZZ/2$ and
$H_8(\BPL,\BO)\cong\ZZ/28$.

\proof The corresponding homotopy groups are given in
Kirby-Siebenmann \cite{KS} p.~246.
\qed
\end{Lem}
From the universal coefficient theorem, see Spanier
\cite{Span} Ch.~5.5 Theorem 3, we have the following
result.
\begin{Prop}
For any coefficient domain $R$, there are exact sequences
\begin{diagram}
H^4(\BO;R)&\lTo& H^4(\BTOP;R)&\lTo^{\tau_{\TOP/\O}}& H^3(\TOP/\O;R)&\lTo & 0\\
H^4(\BPL;R)&\lTo& H^4(\BTOP;R)&\lTo^{\tau_{\TOP/\PL}}& H^3(\TOP/\PL;R)&\lTo & 0\\
H^8(\BO;R)&\lTo& H^8(\BPL;R)&\lTo^{\tau_{\PL/\O}}& H^7(\PL/\O;R)&\lTo & 0.\\
\end{diagram}
\qed
\end{Prop}
\begin{Num}
\label{ExoticClasses}
The $\tau_{\TOP/\PL}$-image
of the generator of $H^3(\TOP/\PL;\ZZ/2)\cong\ZZ/2$ is the universal
\emph{Kirby-Siebenmann class} $ks\in H^4(\BTOP;\ZZ/2)$. For $R=\ZZ/4$,
we pick a generator $\kappa$ of the $\tau_{\PL/\O}$-image
of $H^7(\PL/\O;\ZZ/4)\cong\ZZ/4$. Thus we define the first exotic
characteristic classes
\[
ks\in H^4(\BTOP;\ZZ/2)\quad\text{ and }\quad
\kappa\in H^8(\BPL;\ZZ/4)
\]
The map $\tau$ is the transgression (cp. McCleary \cite{MC} 6.2).
This is maybe not obvious
from our construction. To see this, consider the following diagram
of cochain complexes.
{\small\begin{diagram}[1.5em]
&& 0 && 0 && 0 \\
&&\uTo && \uTo && \uTo \\
0 & \lTo & S^\bullet(\PL/\O) & \lTo &
S^\bullet(\BO) & \lTo & S^\bullet(\BO,\PL/\O) & \lTo & 0 \\
&&\uTo && \uTo && \uTo \\
0 & \lTo & S^\bullet(C_{\PL/\O}) & \lTo &
S^\bullet(\BPL) & \lTo & S^\bullet(\BPL,C_{\PL/\O}) & \lTo & 0 \\
&&\uTo && \uTo && \uTo \\
0 & \lTo & S^\bullet(C_{\PL/\O},\PL/\O) & \lTo &
S^\bullet(\BPL,\BO) & \lTo & S^\bullet(\BPL,\BO,C_{\PL/\O}) & \lTo & 0 \\
&&\uTo && \uTo && \uTo \\
&& 0 && 0 && 0 
\end{diagram}}%
Here, $S^\bullet(\BPL,\BO,C_{\PL/\O})$ denotes the singular cochain complex
of the triad $(\BPL,\BO,\linebreak C_{\PL/\O})$, see Eilenberg-Steenrod 
\cite{ES} VII.11.
Using the isomorphisms derived above,
patient diagram chasing in the corresponding big
diagram for cohomology (the infinite cohomology jail window, see 
eg. Cartan-Eilenberg \cite{CE} IV Proposition 2.1) shows
that $\tau$ is indeed the transgression.
\end{Num}

\section{Stable $\RR^n$-bundles over $\SS^m$}
\label{BundleClass}
Our aim is the classification of $\RR^m$-bundles over
$\SS^m$ in terms of characteristic classes.
We begin with the stable classification, which is easier.
Recall from \ref{HirschExact} that
there is an exact sequence
\[
0\rTo \pi_{k}(\BO)\rTo \pi_{k}(\BTOP)\rTo\pi_{k-1}(\TOP/\O)\rTo 0
\]
for all $k\geq 0$.
\begin{Lem}
In dimensions $k=2,4,8$ these exact sequences read as follows.
\begin{diagram}
0&\rTo &\ZZ/2&\rTo&\ZZ/2&\rTo&0&\rTo&0 && (k=2)\\
0&\rTo &\ZZ&\rTo&\ZZ\oplus\ZZ/2&\rTo&\ZZ/2&\rTo&0&& (k=4)\\
0&\rTo &\ZZ&\rTo&\ZZ\oplus\ZZ/4&\rTo&\ZZ/28&\rTo&0&& (k=8)
\end{diagram}

\proof
We consider the maps
\begin{diagram}[size=2em]
& &  & & \G/\TOP \\
      &      &  &      & \dTo \\
\TOP/\O & \rTo & \BO & \rTo & \BTOP \\
          &      &  &     & \dTo \\
          &      & && \BG 
\end{diagram}
The space $\TOP/\O$ is $2$-connected, see
Kirby-Siebenmann \cite{KS} p.~246,
and this establishes the result for $m=2$. Furthermore
$\pi_3(\TOP/\O)\cong\ZZ/2$ and $\pi_7(\TOP/\O)\cong\ZZ/28$, see 
Kirby-Siebenmann \cite{KS} p.~246 and 200, and  Kervaire-Milnor
\cite{KM}.
From the isomorphisms $\pi_k(\BG)\cong\pi_{k-1}^s(\SS^0)$
we have $\pi_4(\BG)\cong\ZZ/24$ and $\pi_8(\BG)\cong\ZZ/240$,
see Toda \cite{Toda} Ch.~XIV, Hu \cite{Hu} Ch.~XI. Theorem~16.4 and p.~332, or
Fomenko-Fuchs-Gutenmacher \cite{FFG} p.~300.
Finally, $\pi_{4k}(\BO)\cong\ZZ$ for $k\geq 1$ by Bott periodicity.
Thus we obtain diagrams
\begin{diagram}
     &      & \ZZ && 
&&   &      & \ZZ \\
     &      & \dTo_{\text{mono}} && 
&&   &      & \dTo_{\text{mono}} \\
\ZZ & \rTo^{\text{mono}} & \pi_4(\BTOP) & \rTo^{\text{epi}} & \ZZ/2 &&
\ZZ & \rTo^{\text{mono}} & \pi_8(\BTOP) & \rTo^{\text{epi}} & \ZZ/28 \\
& \rdTo_{J_\O}^{\text{epi}}  & \dTo_{\text{epi}} &
&&&& \rdTo_{J_\O}^{\text{epi}}  & \dTo_{\text{epi}} \\
 & & \ZZ/24 &&&&
 & & \ZZ/240. \\
\end{diagram}
In these diagrams, the rows are short exact sequences by the remarks
at the beginning of this section. The columns are also short exact,
since $\pi_{4k-1}(\G/\TOP)=0$, while $\pi_{4k+1}(\BG)$ is finite.
The slanted arrows $J_\O$ are known to be
epimorphisms in these dimensions,
see eg.~Adams \cite{AdamsJ-IV} p.~22 and p.~46.

Let $T_{4k}$ denote the torsion group of $\pi_{4k}(\BTOP)$, i.e
$\pi_{4k}(\BTOP)\isom\ZZ\oplus T_{4k}$.
Suppose that $m=4$. The diagram shows that $T_4$ injects into $\ZZ/2$.
If we tensor the diagram with $\ZZ/2$, then $J_\O\otimes\ZZ/2$
is a bijection. Therefore, the horizontal sequence
\[
0\too\ZZ/2\too\ZZ/2\oplus(T_4\otimes\ZZ/2)\too\ZZ/2\too0
\]
is still exact after tensoring
(the only point to check is injectivity of the second
arrow). It follows that $T_4\isom\ZZ/2$.

The case $m=8$ is similar. The group $T_8$ injects
into $\ZZ/28$ and into $\ZZ/240$.
Since $\gcd(28,\linebreak 240)=4$, the group
injects into $\ZZ/4$. Similarly as in the case $m=4$, tensoring the
diagram with $\ZZ/4$ we see that $J_\O\otimes\ZZ/4$ is an
isomorphism. Thus the sequence
\[
0\too\ZZ/4\too\ZZ/4\oplus(T_8\otimes\ZZ/4)\too\ZZ/4\too0
\]
is exact, and $T_8\isom\ZZ/4$.
\qed
\end{Lem}
For $k\geq 3$, the structure of the torsion groups $T_{4k}$ was
determined by Brumfiel \cite{Brumfiel}, see Madsen-Milgram \cite{MM} p.~117.
The cases $k=1,2$ are special; they are considered in
Kirby-Siebenmann \cite{KS} p.~318 and Williamson \cite{Wil} p.~29.
\begin{Num}
\label{PontrjaginCoKernel}
The result above yields thus exact sequences
\begin{diagram}[width=3.5em]
0& \rTo & \pi_{4}(\BO) & \rTo^\cong  & \pi_{4}(\BTOP)/T_4 &
 \rTo & 0 & \rTo & 0 & (m=4) \\
 0& \rTo & \pi_{8}(\BO) & \rTo  & \pi_{8}(\BTOP)/T_8 &
 \rTo & \ZZ/7 & \rTo & 0 & (m=8).
\end{diagram}
The cokernels of the corresponding maps 
$\pi_{4k}(\BO)\too\pi_{4k}(\BTOP)/T_{4k}$ for $k\geq 3$
are determined in Brumfiel \cite{Brumfiel} p.~304 in
number theoretic terms.
\end{Num}
Consider now the map which assigns
to a stable $\RR^n$-bundle $\xi$ over $\SS^{4k}$
the rational Pontrjagin number $\bra{p_{4k}(\xi),[\SS^{4k}]}$.
We want to determine the possible values of this map,
and we do this first for vector bundles.

\begin{Num}
For a finite connected
CW-complex $X$ and an $n$-dimensional vector bundle
$\xi$ over $X$,
there is a classifying map $X\too \BO(n)$. The Pontrjagin classes 
of $X$ are obtained by pulling back the universal Pontrjagin
classes in $H^\bullet(\BO;\QQ)$ via the composite
$X\too\BO(n)\too\BO$. Only the homotopy type of this map
is important, so we view it as an element of the set $[X;\BO]$
of free homotopy classes of maps from $X$ to $\BO$
(since $\BO$ is an $H$-space and $X$ is connected,
we have $[X;\BO]_0=[X,\BO]$, cp.~the remarks in the introduction
of this paper). Now this
set $[X;\BO]$ can be identified with the real
reduced $\widetilde{KO}$-theory of $X$,
\[
[X;\BO]=\widetilde{KO}(X),
\]
see Atiyah-Hirzebruch \cite{AH}, Hirzebruch \cite{Hirz}
or Husemoller \cite{Hus}.
The \emph{Pontrjagin character} is the ring
homomorphism $ph=ch\circ cplx$,
\[
{KO}(X){\too^{cplx}}
{KU}(X){\too^{ch}}
H^\bullet(X;\QQ)
\]
where $cplx$ denotes complexification of real vector bundles, and
$ch$ denotes the \emph{Chern character} of complex
$KU$-theory, see Atiyah-Hirzebruch \cite{AH} or Hirzebruch \cite{Hirz} 1.4.
For the following facts see Hirzebruch \cite{Hirz} 1.4--1.6.
Hirzebruch's integrality theorem says that
\[
ch(KU(\SS^{2k}))=H^\bullet(\SS^{2k})\subseteq H^\bullet(\SS^{2k};\QQ),
\]
with $ch(\eta)=\mathrm{rk}_\CC(\eta)+(-1)^{k-1}\frac1{(k-1)!}c_k(\eta)$.
Recall also that $p_{4k}(\xi)=(-1)^{k}c_{2k}(\xi\otimes\CC)$.
Thus we have the formula
\[
ph(\xi)=\mathrm{rk}_\RR(\xi)+(-1)^{k-1}\frac1{(2k-1)!}p_{4k}(\xi)
\]
on $\SS^{4k}$. The map
$\widetilde{KO}(\SS^{4k})\too^{cplx}\widetilde{KU}(\SS^{4k})$
is injective, with
cokernel $0$ for $k$ even, and cokernel $\ZZ/2$ for $k$ odd.
\end{Num}
Combining these facts, we have the following result.
\begin{Lem}
\label{IntegralityOnSphere}
Let $\xi$ be a vector bundle over $\SS^{4k}$, and let
$x\in H^{4k}(\SS^{4k})$ be a generator. Then
\[
p_{4k}(\xi)=a_\xi\cdot d_k\cdot(2k-1)!\cdot x
\]
where $a_\xi$ is an integer depending on $\xi$, and $d_k=1$ for $k$ even,
$d_k=2$ for $k$ odd. Conversely, given any integer $a_\xi$, there exists
a vector bundle $\xi$ with such a Pontrjagin class, and $\xi$ is
unique up to stable equivalence.

\proof
Since $\BO$ is an $H$-space, the set $[\SS^{4k};\BO]$
of free homotopy classes coincides with the homotopy group
$\pi_{4k}(\BO)=[\SS^{4k};\BO]_0$. Viewing the universal
Pontrjagin class $p_{4k}$ as a map $\BTOP\rTo^{p_{4k}} K(\QQ,4k)$
into an Eilenberg-MacLane space,
we obtain a nontrivial group homomorphism
\[
\ZZ\cong\pi_{4k}(\BO)
\rTo^{(p_{4k})_\#}\pi_{4k}(K(\QQ,4k))\cong\QQ.
\]
\qed
\end{Lem}
If $\xi$ is a vector bundle, 
the possible values of the rational numbers
$\bra{p_{4k}(\xi),[\SS^m]}$ are thus the integral multiples of
$d_k\cdot(2k-1)!$.
\begin{Lem}
\label{PontrjaginClassesLemma}
Let $\xi$ be an $\RR^n$-bundle over $\SS^m$, for $m=4,8$. Then
\begin{align*}
\bra{p_{4}(\xi),[\SS^4]}&\in 2\ZZ\subseteq\QQ \\
\bra{p_{8}(\xi),[\SS^8]}&\in\textstyle \frac67\ZZ\subseteq\QQ .
\end{align*}
Conversely, for each of these values, there exists an $\RR^n$-bundle
(for some sufficiently large $n$)
whose Pontrjagin number assumes this value.

\proof
This is clear from Lemma \ref{IntegralityOnSphere}, applied to the
special cases $k=1,2$, and the formula for the
cokernel of the map $\pi_{4k}(\BO)\too\pi_{4k}(\BTOP)/T_{4k}$
which was derived in \ref{PontrjaginCoKernel}.
\qed
\end{Lem}
We have proved the following result.
\begin{Prop}
Let $\xi$ be an $\RR^n$-bundle over $\SS^m$, for $m=2,4,8$.
Up to stable equivalence, the bundle $\xi$ is completely determined
by the following characteristic classes:
\begin{description}
\item[($m=2$)] its $2$nd Stiefel-Whitney class $w_2(\xi)$.
\item[($m=4$)] its $4$-dimensional Pontrjagin class $p_4(\xi)$ and its
Kirby-Siebenmann class $ks(\xi)$.
\item[($m=8$)] its $8$-dimensional Pontrjagin class $p_8(\xi)$ and the 
characteristic class $\kappa(\xi)$.
\end{description}
The possible ranges for the values of
these characteristic classes, evaluated on the fundamental class $[\SS^m]$,
are $\ZZ/2$ (for $m=2$), $2\ZZ$ and $\ZZ/2$ (for $m=4$),
and $\frac67\ZZ$ and $\ZZ/4$ (for $m=4$), respectively.

\proof
We prove the $8$-dimensional case; the others are similar.
Let $(\BPL,\BO)\rTo\linebreak (K(\ZZ/4,8),*)$ represent the generator
of $H^8(\BPL,\BO;\ZZ/4)$ which maps to $\kappa$.
Then the composite $\BPL\rTo(\BPL,\BO)\rTo (K(\ZZ/4,8),*)$
induces an isomorphisms on the $\ZZ/4$-factors
in $\pi_8(\BPL)\rTo\pi_8(\BPL,\BO)\rTo\pi_8(K(\ZZ/4,8))$, which
can be identified with the map $\xi\mapstoo \kappa(\xi)\in H^8(\BPL;\ZZ/4)$.
The octonionic Hopf bundle $\eta_\OO$
over $\SS^8$ represents an element
of $\pi_8(\BPL)$ with $p_8(\eta_\OO)\neq 0$. Thus, the map
$\xi\mapsto (p_8(\xi),\kappa(\xi))\in H^8(\SS^8;\QQ)\oplus H^8(\SS^8;\ZZ/4)
\cong\QQ\oplus\ZZ/4$ is an injection (with image $\frac67\ZZ\oplus\ZZ/4$).
\qed
\end{Prop}

\section{$\RR^m$-bundles over $\SS^m$}
In the previous section, we classified bundles over
$\SS^m$ up to stable equivalence in terms of characteristic
classes.
To obtain an unstable classification, i.e. a classification of
$\RR^m$-bundles over $\SS^m$, we use the following result.
\begin{Prop}
\label{KernelProp}
Let $m\geq 2$ be even. Then there is a commutative diagram
with exact rows
\begin{diagram}[height=2em]
0 & \rTo & \pi_m(\SS^m) & \rTo & \pi_m(\BO(m)) & \rTo &
\pi_m(\BO) & \rTo & 0 \\
&& \dEq && \dTo && \dTo \\
0 & \rTo & \pi_m(\SS^m) & \rTo & \pi_m(\BPL(m)) & \rTo &
\pi_m(\BPL) & \rTo & 0 \\
&& \dEq && \dTo && \dTo \\
0 & \rTo & \pi_m(\SS^m) & \rTo & \pi_m(\BTOP(m)) & \rTo &
\pi_m(\BTOP) & \rTo & 0 \\
&& \dEq && \dTo && \dTo \\
0 & \rTo & \pi_m(\SS^m) & \rTo & \pi_m(\BG(m)) & \rTo &
\pi_m(\BG) & \rTo & 0 .
\end{diagram}
In this diagram, the second column of vertical arrows is
induced by the respective classifying map of the tangent
bundle of $\SS^m$ (resp. its underlying spherical fibration)
and the third column of vertical arrows corresponds to
stabilization.

\proof
This is proved in \cite{KrP=L} pp.~93--95; there, the result
is stated for the oriented case, but the homotopy groups
are the same. The PL result is not stated in \cite{KrP=L},
but in low dimensions, $\pi_k(\BO)\cong\pi_k(\BPL)$ and
in higher dimensions $\pi_k(\BPL)\cong\pi_k(\BTOP)$,
see Kirby-Siebenmann \cite{KS} V \S5. For a related result 
(but excluding dimension $4$) see Varadarajan \cite{Var}.
\qed
\end{Prop}
Our aim is to prove that the first three rows in this diagram
split, provided that $m\geq 4$. It suffices to
prove this for the third row; the diagram then implies the splitting
of the first two rows.
\begin{Lem}
\label{SplittingLemma}
Let $m\geq 4$ be even. Then the exact sequence
\[
0 \rTo \pi_m(\SS^m) \rTo \pi_m(\BTOP(m)) \rTo 
\pi_m(\BTOP)  \rTo  0 
\]
splits.

\proof
Since we consider only higher dimensional homotopy groups,
we may as well consider the universal coverings $\BSTOP(m)$ and $\BSTOP$
which classify oriented bundles. So we have an exact sequence
\[
0\rTo\pi_m(\SS^m)\rTo\pi_m(\BSTOP(m))\rTo\pi_m(\BSTOP)\rTo0.
\]
Let $e$ denote the universal Euler class, viewed as a map
$\BSTOP(m)\too^e K(\ZZ,m)$ into an Eilenberg-MacLane space.
The Euler class yields thus a homomorphism
\[
\pi_m(\BSTOP(m))\too^{e_\#}\pi_m(K(\ZZ,m))=H^m(\SS^m)
\too^{\bra{-,[\SS^m]}}_\cong\ZZ.
\]

If $m\neq 4,8$, then the image of this map is $2\ZZ$ by
Adams' result, cp.~Proposition \ref{IfManifolde=1}.
The Euler class of the tangent bundle of $\SS^m$ is $2x$,
so the exact sequence above splits: we have constructed a left inverse for
the second arrow $\pi_m(\SS^m)\too\pi_m(\BSTOP(m))$.

For $m=4,8$, we consider the last two rows of the big diagram,
\begin{diagram}[height=2em]
0 & \rTo & \pi_m(\SS^m) & \rTo & \pi_m(\BTOP(m)) & \rTo
& \pi_m(\BTOP) & \rTo & 0 \\
 & & \dEq && \dTo && \dTo \\
0 & \rTo & \pi_m(\SS^m) & \rTo & \pi_m(\BG(m)) & \rTo
& \pi_m(\BG) & \rTo & 0. 
\end{diagram}
Let $T_m'$ denote the torsion subgroup in $\pi_m(\BTOP(m))$,
and $T_m\cong\ZZ/(m/2)$ the torsion subgroup of $\pi_m(\BTOP))$.
It is clear that $T_m'$ injects into $T_m$. The
sequence splits if and only if $T_m'$ maps isomorphically onto
$T_m$.
Now $\pi_4(\BG(4))\isom\pi_3(\G(4))\isom\pi_7(\SS^4)
\isom\ZZ\oplus\ZZ/12$, and $\pi_4(\BG)\isom\pi_3^s(\SS^0)\isom\ZZ/24$.
Similarly,
$\pi_8(\BG(8))\isom\pi_7(\G(8))\isom\pi_{15}(\SS^8)
\isom\ZZ\oplus\ZZ/120$, and $\pi_8(\BG)\isom\pi_7^s(\SS^0)\isom\ZZ/240$.
See Toda \cite{Toda} Ch.~XIV, Hu \cite{Hu} Ch.~XI. Theorem~16.4 and p.~332,
or Fomenko-Fuchs-Gutenmacher \cite{FFG} p.~300 for these groups.
It follows that in the bottom row of the diagram,
a generator $\iota_m\in\pi_m(\SS^m)$
maps to $(2,-1)\in\ZZ\oplus\ZZ/12$ resp. in $\ZZ\oplus\ZZ/120$, for $m=4,8$.
If we tensor the diagram above with $\ZZ/2$ for $m=4$ (resp. with
$\ZZ/4$ for $m=8$), then the image of $\iota_m$ still has order 2
(resp. 4) in $\pi_m(\BG(m))\otimes\ZZ/(m/2)$. Thus the bottom row remains
exact after
tensoring, and therefore, the upper row remains also exact (the only
point to be checked was the injectivity of the second arrow).
It follows that $T_m'\cong T_m$.
\qed
\end{Lem}
In dimension $m=2$, we use Kneser's old result $\BO(2)\simeq\BTOP(2)$,
see Kirby-Siebenmann \cite{KS} p.~254.
Thus, there is a (non-split) short exact sequence
\[
0\rTo\pi_2(\SS^2)\rTo\pi_2(\BTOP(2))\rTo\ZZ/2\rTo0
\]
and $\pi_2(\BTOP(2))\cong\ZZ$.
Using similar ideas as above, it is not difficult to prove
that the sequence
\[
0\rTo\pi_m(\SS^m)\rTo\pi_m(\BG(m))\rTo\pi_m(\BG)\rTo0
\]
splits for all even $m\neq 2,4,8$ (and for these three values, the
sequence is not split). We will not need this result.

\begin{Num}
\label{StabilizingPropProof}
\emph{Proof of Proposition \ref{StabilizingProp}}
Let $\xi$ and $\xi'$ be $\RR^m$-bundles over $\SS^m$, for
$m\geq 2$ even. Assume first that we can choose orientations
of these bundles such that $e(\xi)=e(\xi')$. Let $c$ denote
the classifying map for the oriented
tangent bundle $\tau\SS^m$, and
let $c_\xi$ and $c_{\xi'}$ be classifying maps for
the oriented bundles $\xi$ and $\xi'$. The splitting
of the exact sequence
\[
0\rTo\pi_m(\SS^m)\rTo^{c_\#}\pi_m(\BSTOP(m))\rTo\pi_m(\BSTOP)\rTo0
\]
implies then that $(c_\xi)_\#=(c_{\xi'})_\#$ (as maps
$\pi_m(\SS^m)\rTo\pi_m(\BSTOP)$). Thus $\xi\cong\xi'$.

In the general case $m\geq 4$, we have the action of
$\pi_1(\BO(m))\cong\pi_1(\BTOP(m))\cong\ZZ/2$ on the higher homotopy groups.
The generator $\alpha_0$ of the fundamental group of $\BO(m)$ maps
$c_\#$ to its negative $-c_\#$, cp.~Steenrod \cite{Steenrod} 23.11.
From the splitting of the exact sequence
\[
0\rTo\pi_m(\SS^m)\rTo\pi_m(\BTOP(m))\rTo\pi_m(\BTOP)\rTo0
\]
and the diagram in Proposition \ref{KernelProp} we see that
$\alpha_0$ changes the sign of the Euler class. Thus, if
$|e|=|e'|$, then we may as well assume that $e(\xi)=e(\xi')$.

The case $m=2$ follows directly from $\BO(2)\simeq \BTOP(2)$.
\qed
\end{Num}
We summarize our classification of $\RR^m$-bundles over $\SS^m$, for
$m=2,4,8$, as follows.
\begin{Prop}
\label{TheBundlesWeObtain}
Let $\xi$ be an $\RR^2$-bundle over $\SS^2$. Up to equivalence,
$\xi$ is determined by its absolute Euler number
$|e|=|\bra{e(\xi),[\SS^2]}|$, and
for each $|e|\in\mathbb N$, there exists one such bundle.
A weak equivalence between any two such bundles is an equivalence.

Let $\xi$ be an $\RR^4$-bundle over $\SS^4$.
Up to equivalence,
$\xi$ is determined by its absolute Euler number
$|e|=|\bra{e(\xi),[\SS^2]}|$, its
Kirby-Siebenmann number $\bra{ks(\xi),[\SS^4]}\in\ZZ/2$
and the Pontrjagin number $\bra{p_4(\xi),[\SS^4]}\in 2\cdot\ZZ$.
For each triple $(|e|,ks,p_4)\in\mathbb N\times\ZZ/2\times 2\cdot\ZZ$
satisfying the relation $p_4+2|e|\equiv0\pmod 4$,
there exists one such bundle.
If two such bundles $\xi,\xi'$ are weakly equivalent, but
not equivalent, then $(|e|,p_4,ks)=(|e'|,-p_4',ks')$.

Let $\xi$ be an $\RR^8$-bundle over $\SS^8$.
Up to equivalence,
$\xi$ is determined by its absolute Euler number
$|e|=|\bra{e(\xi),[\SS^2]}|$, the
number $\bra{\kappa(\xi),[\SS^8]}\in\ZZ/4$
and the Pontrjagin number $\bra{p_8(\xi),[\SS^8]}\in \frac67\cdot\ZZ$.
For each triple in $\mathbb N\times\ZZ/4\times\frac67\cdot\ZZ$
satisfying the relation $\frac73p_8+2|e|\equiv0\pmod4$,
there exists one such bundle.
If two such bundles $\xi,\xi'$ are weakly equivalent, but
not equivalent, then $(|e|,p_8,\kappa)=(|e'|,-p_8',-\kappa)$.

\proof
We prove the $8$-dimensional case; the other cases are similar.
First, we classify oriented bundles; the orientation we choose is
the orientation determined by the universal oriented $\RR^8$-bundle
over $\BSTOP(8)$. Then it is clear from our discussion that
$\xi$ is determined by the data
\[
\textstyle
(\bra{e(\xi),[\SS^8]},\bra{\kappa(\xi),[\SS^8]},\bra{p_8(\xi),[\SS^8]})
\in\ZZ\times\ZZ/4\times \frac67\ZZ.
\]
Now we have as in \ref{StabilizingPropProof} the action of $\alpha_0$
which changes the sign of the Euler class without changing the sign
of $p_8$ and $\kappa$ (since these two classes come from $\BTOP$,
where $\alpha_0$ acts trivially). This shows that the
given numbers classify the bundle  up to equivalence.

Let $\iota_8\in\pi_8(\SS^8)$ denote the canonical generator, and let
$c$ be the classifying map for the oriented tangent bundle of
$\SS^8$.
For the number-theoretic relation between the Pontrjagin class
and the Euler class, we note first that the image
$(c)_\#(\pi_8(\SS^8))$ is a direct factor in $\pi_8(\BSTOP(8))$.
The octonionic
Hopf line bundle $\eta_\OO$ has Euler class $x$ and
Pontrjagin class $6x$ (for a suitable generator $x$ of $H^m(\SS^m)$),
see eg. \cite{KrSmooth} Theorem 9.
Let $h$ be a classifying map for the
oriented bundle $\eta_\OO$.
Then $c_\#(\iota_8)$ and $h_\#(\iota_8)$
span $\pi_8(\BSTOP(8))\otimes\QQ$ (over $\QQ$), and
the image of $h$ in $\pi_8(\BSTOP)$ spans $\pi_8(\BSTOP)\otimes\QQ$
(over $\QQ$). Since $e(\xi)$ is necessarily integral for any
oriented $\RR^8$-bundle over $\SS^8$, we see from
Lemma \ref{PontrjaginClassesLemma} that there
exists a bundle $\eta'$ with classifying map $h'$ and
$e(\eta')=x$, $p_8(\eta')=\frac67 x$, and that
$c_\#(\iota_8)$ and $h'_\#(\iota_8)$ span a direct complement
of the torsion group of $\pi_8(\BSTOP(8))$ (over $\ZZ$).

Finally,
a weak equivalence which is not an equivalence comes from a map
$\SS^8\rTo\SS^8$ of degree $-1$; such a map changes the sign of
every characteristic class.
\qed
\end{Prop}

\section{The classification of the models}
\label{ClassifyModels}
In this section, we obtain the final homeomorphism classification
of our models. We fix some notation.
Let $\xi$ be an $\RR^m$-bundle over $\SS^m$, with
absolute Euler number $|e|=1$, for $m=2,4,8$,
let $M(\xi)$ be its Thom space,
and let $s_0:\SS^m\rTo E\subseteq M(\xi)$ be the zero-section.
Let $y_m\in H^m(M(\xi))$ denote a generator, such that
$x=s_0^\bullet y_m$ is a generator dual to the chosen orientation
$[\SS^m]$. For $m=4,8$, we have by \ref{StableClassInHalfDim} the relations
\begin{align*}
s_0^\bullet(p_m(M(\xi))) & =p_m(\xi)      &&\quad (m=4,8)\\
s_0^\bullet(ks(M(\xi)))&=ks(\xi)          &&\quad (m=4)\\
s_0^\bullet(\kappa(M(\xi)))&=\kappa(\xi)  &&\quad (m=8).
\end{align*}
\begin{Thm}
\label{TheModelsWeObtain}
Up to homeomorphism,
our construction yields precisely the following models.

For $m=2$, there is just one model, the complex projective
plane $M(\eta_\CC)\cong\CP^2$, where $\eta_\CC$ is the complex Hopf
line bundle (the tautological bundle) over $\CP^1=\SS^2$.

For $m=4$, let $p_4(\xi)=2(1+2t)x$ and $ks(\xi)=sx$, for
$(t,s)\in\ZZ\times\ZZ/2$. By Proposition \ref{TheBundlesWeObtain},
the pair $(r,s)=(1+2t,s)$
determines $\xi$ up to equivalence, so we may put $M_{r,s}=M(\xi)$.
If there is a homeomorphism $M_{r,s}\cong M_{r',s'}$, then
$(r,s)=(\pm r',s)$. The quaternion projective plane
is $M_{1,0}\cong\HP^2=M(\eta_\HH)\cong M_{-1,0}$. The model $M_{r,s}$ admits
a PL structure (unique up to isotopy) if and only if $s=0$.

For $m=8$, let $p_8(\xi)=\frac67(1+2t)x$ and $\kappa(\xi)=sx$, for
$(t,s)\in\ZZ\times\ZZ/4$. Again, the pair $(r,s)=(1+2t,s)$ determines
$\xi$ up to equivalence by Proposition \ref{TheBundlesWeObtain},
so we may put $M_{r,s}=M(\xi)$.
If there is a homeomorphism $M_{r,s}\cong M_{r',s'}$, then
$(r,s)=\pm(r',s')$. The octonionic projective plane
is $M_{7,0}\cong\OP^2=M(\eta_\OO)\cong M_{-7,0}$. Each model $M_{r,s}$ admits
a PL structure (unique up to isotopy).

\proof
The topological result is clear from our classification of
$\RR^m$-bundles over $\SS^m$ with absolute Euler number $|e|=1$.
For a closed manifold $M$ of dimension at least $5$, the only obstruction
to the existence of a PL structure is the Kirby-Siebenmann class $ks(M)$.
If such a PL structure exists on $M$ , the number of isotopy classes
of PL structures is determined by $[M;\TOP/\PL]\cong H^3(M;\ZZ/2)$,
see Kirby-Siebenmann \cite{KS} p.~318.
In our case, these groups are zero.
\qed
\end{Thm}
The models constructed by Eells-Kuiper in \cite{EK}
are $\CP^2$ (for $m=2$),
the models $M_{r,0}$, with $r=1+2t$ (for $m=4$), and the models
$M_{7r,0}$, with $r=1+2t$ (for $8=4$). Brehm-K\"uhnel \cite{BK}
construct $8$-dimensional PL manifolds which look like projective planes,
and with small numbers of vertices. However, the
characteristic classes of their examples seem to be unknown.
\begin{Num}
Now we consider the question which of the models admit a DIFF
structure.
The number of PL or DIFF structures on $\CP^2$ is presently not
known, so we concentrate from now on the cases $m=4,8$.
Concerning DIFF structures in higher dimensions,
a necessary condition (besides the existence
of a PL structure) is clearly
that $\tau M(\xi)$ admits a vector bundle
structure. In particular, $\xi\cong s_0^*\tau M(\xi)$ has to admit
a vector bundle structure. Thus, if $M_{r,s}=M_{1+2t,s}$
admits a DIFF structure, then clearly $s=0$, and in addition
$t\equiv 3\pmod 7$ for $m=8$. But this guarantees only that
$\xi$ admits a vector bundle structure.
From \ref{SpecificValues} and Lemma \ref{RelationsForModelsLem},
we see that
\begin{align*}
\Ahat\ [M_{1+2t,s}]&=-\frac{t(1+t)}{2^3\cdot 7}&&(m=4)\\
\Ahat\ [M_{7(1+2u),s}]&=-\frac{u(1+u)}{2^7\cdot 127}&&(m=8)
\end{align*}
Since $M_{r,s}$ is $3$-connected, the first Stiefel-Whitney classes
vanish. Thus, if $M_{r,s}$ admits a DIFF structure, then it is
a Spin manifold, see Lawson-Michelsohn \cite{LM} Ch.~II Theorem 2.1.
But for a closed oriented Spin manifold $M^{4k}$, the $\Ahat$-genus
$\Ahat[M]$ is precisely the index of the Atiyah-Singer operator,
see Lawson-Michelsohn \cite{LM} Ch.~IV Theorem 1.1;
in particular, it is an integer.
\end{Num}
\begin{Lem}
If $M_{r,s}$ admits a DIFF structure, then
we have the following relations:
For $m=4$ put $r=1+2t$. Then $s=0$ and $t\equiv 0,7,48,55\pmod{56}$.
For $m=8$ put $r=7(1+2u)$. Then $s=0$ and $u$ is an integer with
$u\equiv 0,127,16128,16255\pmod{16256}$.
\qed
\end{Lem}
This corresponds to Eells-Kuiper \cite{EK} Proposition 10 A, p.~43.
In fact, the above conditions are sharp.
\begin{Thm}
If $m=4$, then $M_{1+2t,s}$ admits a DIFF structure if and
only if $s=0$ and $t\equiv 0,7,48,55\pmod {56}$.
If $m=8$, then $M_{7(1+2u),s}$ admits a DIFF structure if and only if
$s=0$ and $u\equiv 0,127,16128,16255\pmod{16256}$.

\proof
The number theoretic condition guarantees that $\xi$ admits a
vector bundle structure. Thus we may choose a Riemannian metric
on $\xi$.
Let $S\!E\too\SS^m$ denote the corresponding unit sphere bundle of $\xi$.
Then $E_0\simeq S\!E$ is a homotopy $2m-1$-sphere and thus
homeomorphic to $\SS^{2m-1}$ by the proof of the generalized
Poincar\'e conjecture, see Smale \cite{Sma} and Newman \cite{New}.
If $S\!E$ is diffeomorphic to $\SS^{2m-1}$,
then we may choose a diffeomorphism $\alpha:\SS^{2m-1}\too S\!E$.
Gluing the closed $2m$-disk $\mathbb D^{2m}$ along $\alpha$
to the closed unit disk bundle $D\!E$ of the vector bundle $\xi$,
we obtain a smooth $2m$-manifold $D\!E\cup_\alpha\mathbb D^{2m}$
homeomorphic to $M_{r,s}\cong D\!E/S\!E$.
Thus the problem is reduced to the question whether $S\!E$ is
diffeomorphic to $\SS^{2m}$. Now the unit disk bundle $X=D\!E$
is an \emph{almost closed manifold}, i.e. a smooth compact manifold $X$ with
boundary $\partial X=S\!E$ a homotopy sphere, see Wall \cite{Walln-1}.
By \emph{loc.cit.} p.~178, such a
manifold $X^{2m}$ has a standard sphere $\SS^{2m-1}$ as its boundary if
and only if its $\Ahat$-genus is integral, provided that $m=4,8$ and
that $X$ is $m-1$-connected. The $\Ahat$-genus of $D\!E$ coincides
of course with the $\Ahat$-genus of $M_{r,s}$.
\qed
\end{Thm}
The case $m=8$ remained open in Eells-Kuiper \cite{EK}.
Note that the existence of a positive scalar curvature metric
implies that the $\Ahat$-genus vanishes, see
Lawson-Michelsohn \cite{LM} Ch.~IV Theorem 4.1;
this happens only for the models
$\HP^2$ and $\OP^2$.
\begin{Prop}
The only models which admit a DIFF structure with a positive scalar
curvature metric are $\HP^2$ and $\OP^2$.
\qed
\end{Prop}
These two manifolds admit a positive scalar curvature metric
for \emph{any} DIFF structure; this follows from Stolz' proof of
the Gromov-Lawson-Rosenberg conjecture, see Stolz \cite{Stolzpscm}
Theorem A (for these two manifolds, Rosenberg's earlier result 
\cite{Ros} for lower dimensions actually suffices). For spin
manifolds of dimension $4k$, the map
$\alpha:\Omega_{4k}^\Spin\rTo KO(\SS^{4k})$ can
be identified with a scalar multiple of the $\Ahat$-genus.

Concerning the number of DIFF structures on a smoothable model
$M_{r,s}$, we have Wall's result \cite{Walln-1}
which says that the DIFF structure on the almost closed manifold
$D\!E$ is unique. The group $\Theta_{2m}$ then acts transitively
on the collection of all smoothings of $M_{r,s}$, see the introduction in
Stolz \cite{Stolz}. This group is cyclic
of order $2$ for $m=4,8$, see Kervaire-Milnor \cite{KM} p.~504.
\begin{Prop}
If a model $M_{r,s}$, for $m=4,8$,
admits a DIFF structure, then it admits at most
two distinct DIFF structures.
\qed
\end{Prop}
This fact was already observed in Eells-Kuiper \cite{EK}.
To obtain a more precise result, i.e. the exact number of
DIFF structures, one would have to determine the
inertia groups of the smoothable models, cp. Stolz \cite{Stolz}.
Maybe his high-dimensional techniques can be adapted to this
situation.

\begin{Num}
Finally, we consider oriented bordisms between distinct models.
So suppose that
\[
\partial W=M_{r,s}\cup -M_{r',s'}
\]
is a compact oriented
(topological) bordism between two models.
Clearly, the numbers $p_m^2[M_{r,s}]$,
$ks^2[M_{r,s}]$ (for $m=4$), and $\kappa^2[M_{r,s}]$ (for $m=8$)
are bordism invariants. Thus, the existence of an oriented bordism
implies that $(r^2,s^2)=({r'}^2,{s'}^2)$.
This is good enough to settle
the case $m=4$; here, we conclude that $(r,s)=\pm(r',s')$ and thus
$M_{r,s}\cong M_{r',s'}$, because
$s\in\ZZ/2$ has no sign. For $m=8$ this is not good enough to
conclude that $(r,s)=\pm(r',s')$ because $s$ is $\ZZ/4$-valued.
So we use the standard
fact from topological surgery theory (as developed in Kirby-Siebenmann
\cite{KS}) that such a bordism $W$ can be made $7$-connected.
The classifying map $W\rTo \BSTOP$ for the oriented tangent bundle lifts
thus to the $7$-connected cover $\BSTOP\bra8$. Put
$\pi=\pi_8(\BSTOP)\cong\ZZ\oplus\ZZ/4$ and let $\BSTOP\bra8\rTo K(\pi,8)$
denote the corresponding characteristic map. Let $x\in H^8(K(\pi,8))\cong
\ZZ\oplus\ZZ/4$ be a generator for a free cyclic factor, and let
$q_8$ denote its image in $H^8(\BSTOP\bra8)$. If $X$ is any $7$-connected
CW-complex, and if $\xi$ is a stable $\RR^n$-bundle over $X$, then the
classifying map $X\rTo^c\BSTOP$ lifts
\begin{diagram}[height=2em,width=4em]
&&\BSTOP\bra8 & \rTo^{q_8} & K(\pi,8)\\
&\ruDotsto &\dTo\\
X & \rTo^c & \BSTOP
\end{diagram}
and the class $q_8(\xi)$ is defined.
From the coefficient pairing $\ZZ\otimes\ZZ/4\rTo\ZZ/4$, 
we obtain for any $\RR^n$-bundle over a $7$-connected space
$X$ a $\ZZ/4$-valued $16$-dimensional characteristic class
$q_8\kappa$ and clearly, this class is a bordism invariant
for $7$-connected oriented bordisms.
For $X=\SS^8$, we know that
$q_8(\xi)=\frac76p_8(\xi)$. For our $16$-dimensional
models, we have thus $q_8(M_{r,s})=\frac76 p_8(M)$ by
\ref{StableClassInHalfDim}, and this
is an odd integral multiple of the generator
$y_8\in H^8(M_{r,s})$.
The $7$-connected bordism yields now the additional relation
$\frac76 p_8\kappa[M_{r,s}]=\frac76 p_8\kappa[M_{r',s'}]$,
which implies that $(r,s)=\pm(r',s')$.
\end{Num}
\begin{Prop}
Non-homeomorphic models $M_{r,s}$, $M_{r',s'}$ fall into different
oriented bordism clas\-ses in $\Omega_{2m}^\STOP$.
\qed
\end{Prop}

\section{The homotopy types of the models}
\label{HomotopyOfModels}
For the homotopy classification of our models we use the
\emph{Spivak fibration}. We recall the construction and
refer to Klein \cite{Klein} for more details.
Let $M$ be a closed $1$-connected manifold (the case where $\pi_1(M)\neq 1$
is more involved and will not be important to us).
There exists an embedding $M\rInto \SS^N$, for some sufficiently
large $N$, such that $M$ hat a normal bundle $\nu M$ in $\SS^N$.
In the group $\widetilde{\mathrm{KTOP}}(M)=[M;\BTOP]$,
the bundle $\nu M$ is just
the inverse of the stable bundle class determined by the tangent
bundle $\tau M$, since
$\tau M\oplus \nu M\cong\tau\SS^N|_M=\underline\RR^N$. This
shows that the normal bundle
$\nu M$ is unique up to stable equivalence.

There is a natural map $\SS^N\too^\alpha M(\nu M)$, where $M(\nu M)$
is the Thom space of the normal bundle. Let $u(\nu M)$ be
an orientation class. One shows that for the fundamental classes
of $\SS^N$ and $M$, one has the relation
\[
u(\nu M)\frown \alpha_\bullet[\SS^N]=[M]
\]
(for the right choice of $u(\nu M)$).
Now a result by Spivak \cite{Spivak}
shows that the stable fiber homotopy type of
the underlying spherical fibration $\sigma M$
of the bundle $\nu M$ depends only on the homotopy type of $M$.
\begin{Num}
Let $c:M\too\BTOP$ be a stable classifying map for $\nu M$, and
let $d:M\too\BTOP$ be a stable classifying map for $\tau M$.
Then $c$ is an inverse of $d$ in the abelian group
$\KTOPtilde(M)=[M;\BTOP]$, and so the composites
$M\too^c\BTOP\too\BG$ and $M\too^d\BTOP\too\BG$ are inverse
to each other in the abelian group $[M;\BG]$.
But $M\too^c\BTOP\too\BG$ is a classifying map for $\sigma M$
and depends thus by Spivak's result only on
the homotopy type of $M$. It follows that
the composite $M\too^d\BTOP\too\BG$ is also a homotopy invariant
of $M$.
\end{Num}
\begin{Prop}
If there is a homotopy equivalence $f:M(\xi)\too^\simeq M(\xi')$
between two models, then there is a fiber homotopy equivalence
$\tau M(\xi)\oplus\underline\RR^k\simeq f^*\tau M(\xi')\oplus\underline\RR^k$,
for some $k\geq 0$.
\qed
\end{Prop}
\begin{Num}
Since $s_0:\SS^m\too M(\xi)$ represents a generator of
$\pi_m(M(\xi))$, this has the consequence that there is
a fiber homotopy equivalence
\[
\xi\oplus\underline\RR^k\simeq g^* 
\xi'\oplus\underline\RR^k,
\]
for some homeomorphism $g:\SS^m\too\SS^m$ of degree $\pm1$.
Since $\xi$ has absolute Euler number $|e|=1$, the $m$th Stiefel-Whitney
class of $\xi$ is nontrivial: the $m$th Stiefel-Whitney class
is the mod $2$ reduction of the Euler class,
see Milnor-Stasheff \cite{MS} Proposition 9.5, so
$0\neq w_m(\xi)=x\in H^m(\SS;\ZZ/2)$.
Also, the Stiefel-Whitney class depends only on the stable type of
the spherical fibration of $\xi$.
Let $R_m\subseteq \pi_m(\BG)$ denote the
set of all elements which represent a spherical fibration over
$\SS^m$ with nontrivial $m$th Stiefel-Whitney class. This is a
coset of a subgroup of index $2$ in $\pi_m(\BG)$ (namely the
kernel of the map $\pi_m(\BG)\rTo^{(w_m)_\#}\pi_m(K(m,\ZZ/2))$).
Precomposing the classifying map with a map of degree $-1$,
we achieve a change of sign for all elements in $\pi_m(\BG)$.
The group $ \pi_m(\BG)\cong\pi^s_{m-1}(\SS^0)$ is cyclic of order $2$,
$24$, $240$, for $m=2,4,8$, respectively, see
Toda \cite{Toda} Ch.~XIV, Hu \cite{Hu} Ch.~XI. Theorem~16.4 and p.~332,
or Fomenko-Fuchs-Gutenmacher \cite{FFG} p.~300. Thus we see that
there are at least $1,6,60$ distinct homotopy types which are
realized by our models, for $m=2,4,8$.
In the Appendix we prove that these numbers are the precise
numbers of homotopy types of Poincar\'e duality complexes 
(see \ref{PDCDef}) which look like projective planes.
\end{Num}
\begin{Thm}
\label{HomotopyTypesOfModels}
Every $1$-connected Poincar\'e duality complex which looks like a
projective plane is homotopy equivalent to one of our models.
The homotopy type of a model $M_{r,s}$ can be determined as follows.
For $m=2$, there is just one model and one homotopy type, namely
$\CP^2$.

If $m=4$, then $M_{r,s}\simeq M_{r',s'}$ if and only if
$r+12s\equiv \pm(r'+12s')\pmod{24}$.

If $m=8$, then $M_{r,s}\simeq M_{r',s'}$ if and only if
$r+60s\equiv \pm(r'+60s')\pmod{240}$.
\qed
\end{Thm}

\section{Our set of models is complete}

In this section we prove that every manifold which
looks like a projective plane is homeomorphic to one
of our model manifolds $M(\xi)$ -- except for dimension 4,
where one has precisely two such manifolds, the model
manifold $\CP^2=M(\eta_\CC)$ and in addition
the Chern manifold $\mathrm{Ch}^4$ which is not a Thom space,
see Theorem \ref{dim4}. This is covered by Freedman's
classification \cite{Fr} of closed $1$-connected $4$-manifolds.
\begin{Thm}
\label{dim4}
There are precisely two $1$-connected closed $4$-manifolds $M$ with
$H_\bullet(M)\cong\ZZ^3$, the complex projective plane $\CC\mathrm{P}^2$
and the Chern manifold $\mathrm{Ch}^4$.

\proof
This is stated in Freedman \cite{Fr} pp.~370--372.
Such a manifold is represented by the odd integral
symmetric bilinear form $\omega=(1)$ on $H_2(M)\cong\ZZ$ and
its Kirby-Siebenmann number $ks[M]\in\ZZ/2$. 
\qed
\end{Thm}
From now on, we assume that $m\neq 2$. The main result of this
section is the following.
\begin{Thm}
\label{ThmCompleteModels}
Let $M^{2m}$ be a manifold which is like a projective plane.
If $m\neq 2$, then $M$ is homeomorphic to one of our models
$M(\xi)$.
\end{Thm}
The proof requires surgery techniques, so
we recall the relevant notions. More information can be
found in 
Madsen-Milgram \cite{MM} Ch.~2, in
Kirby-Siebenmann \cite{KS} Essay V App.~B, in
Wall \cite{Wallsurgery} Ch.~10, and in particular in
Kreck \cite{Kreck}. The basic fact to keep in mind is that
by the results of Kirby-Siebenmann \cite{KS}, higher dimensional
surgery works well in the topological category.
The case $m=2$ can in principle be handled by similar methods, see
Freedman-Quinn \cite{FQ}.
\begin{Num}
\label{PDCDef}
The spaces we are dealing with
are 1-connected, and this simplifies some points.
Suppose that $X$ is a finite and 1-connected
CW-complex. Assume moreover that there is an
element $[X]\in H_n(X)$ such that the cap product induces an
isomorphism
\[
H^k(X)\too^{\ \frown[X]}_\cong H_{n-k}(X)
\]
for all $k$. Then $X$ satisfies Poincar\'e duality, and
the pair $(X,[X])$ is what is called a \emph{Poincar\'e duality complex}
(of formal dimension $n$).
Every closed, 1-connected and oriented manifold is a Poincar\'e
duality complex (we consider here only the 1-connected case; the
presence of a fundamental group requires the more complicated notion
of a simple homotopy type).
\end{Num}
\begin{Num}
Next, recall that an \emph{$h$-cobordism} $(W;M_1,M_2)$ is a
simply connected compact bordism between (simply connected)
closed manifolds $M_1,M_2$,
\[
\partial W=M_1\dot\cup M_2,
\]
with the property that
the inclusions $M_1\too W \lTo M_2$
are homotopy equivalences; an example is the product
cobordism, $(M\times[0,1],M\times 0,M\times 1)$. In higher
dimensions, this is in fact the only example.
\end{Num}
\begin{Num}{\bf $h$-cobordism Theorem}
\label{h-Cobordism}
{\em
Every $h$-cobordism $(W;M_1,M_2)$ with $\dim(W)\geq 5$ is
a product bordism, $W\cong M_1\times[0,1]$.
In particular, there is a homeomorphism $M_1\cong M_2$.}
\proof
For $\dim W\geq 6$, this is proved in Kirby-Siebenmann \cite{KS},
but unfortunately not stated explicitly as a Theorem; see
\emph{loc. cit.} p.~113 and p.~320. For $\dim(W)\geq 7$, a proof is given by
Okabe \cite{Okabe}.
The case $\dim(W)=5$ is proved in Freedman-Quinn \cite{FQ},
with some remarks on the higher dimensional case.
\qed
\end{Num}
\begin{Num}
Suppose now that $M$ is a closed oriented manifold of dimension at least $5$,
and that $f:M\too X$ is a homotopy equivalence, with $f_\bullet[M]=[X]$.
Then $f$ is called
a homotopy manifold structure on $X$; two such homotopy
manifold structures
$M_1\too^{f_1}X\lTo^{f_2} M_2$ are
called \emph{equivalent} if there exists an $h$-cobordism
$(W;M_1,M_2)$ and a map $F:W\too X$ such that the diagram
\begin{diagram}
M_1 & \rInto & W & \lInto & M_2 \\
&  \rdTo_{f_1}&\dTo^F& \ldTo^{f_2}\\
&& X
\end{diagram}
commutes. This relation is transitive and symmetric;
the set of all equivalence classes of homotopy manifold structures
on $X$ is the \emph{structure set} $\mathcal{S}_\TOP(X)$.
Since we are assuming that $\dim(M)\geq 5$, the $h$-cobordism Theorem
\ref{h-Cobordism} applies, and thus every element of
$\mathcal{S}_\TOP(X)$
represents a well-defined homeomorphism type of a closed manifold
homotopy equivalent to $X$.

Let $\mathrm{Aut}(X)\subseteq[X;X]$
denote the group of all self-equivalences of $X$. 
If $\mathcal{S}_\TOP(X)$ is nonempty, there is a natural
action of $\mathrm{Aut}(X)$ on $\mathcal{S}_\TOP(X)$, and the
orbit set 
\[
\mathcal{M}_\TOP(X)=\mathcal{S}_\TOP(X)/\mathrm{Aut}(X)
\]
can be identified with the set of all homeomorphism types of manifolds
homotopy equivalent with $X$.
\end{Num}
\begin{Num}
In order to determine the structure set $\mathcal{S}_\TOP(X)$, it
is convenient to introduce yet another set which contains
$\mathcal{S}_\TOP(X)$ as a subset, the set 
$\mathcal{T}_\TOP(X)$ of \emph{tangential invariants}.
Recall from Section \ref{HomotopyOfModels}
that associated to a $1$-connected Poincar\'e duality complex
$X$ is a spherical fibration, the \emph{Spivak normal bundle}
$\sigma X$ whose stable fiber homotopy class
depends only on the homotopy type of $X$. Let $S\tau M$ denote the
spherical fibration of the topological tangent bundle $\tau M$. If
$f:M\too X$ is a homotopy equivalence, then the spherical fibration
$f^*\sigma X\oplus S\tau M$ is stably fiber homotopically
trivial. To put it differently, $f^*\sigma X$ is
the stable inverse of the spherical fibration $S\tau(M)$ 
of the tangent bundle in the fiber homotopy category.
Let $\sigma_T X$ be a stable inverse of
$\sigma X$. Then $f^*\sigma_TX$ is
stably fiber homotopy equivalent to $S\tau(M)$; in particular, the
stable fiber homotopy type of $S\tau(M)$ is a homotopy invariant of $M$.
We used this fact already in Section \ref{HomotopyOfModels}.
To make things more concrete, let us say that $\sigma_TX$
is an $N$-spherical fibration, for $N>2\dim(M)$. Then
$f:M\too X$ induces a bundle map
$\tau M\oplus\underline\RR^{N-\dim(M)}\too f^*\sigma_TX$
which is a fiber homotopy equivalence. Such a map is called a
TOP \emph{reduction} of $f^*\sigma_TX$.
Two reductions are called equivalent if they differ by a fiber homotopy
equivalence; the set of all stable TOP reductions of a spherical fibration
$\phi$ is denoted $\RTOP(\phi)$. One can show that
every element of $\mathcal{S}_\TOP(X)$ yields a well-defined reduction
of $\sigma_T(X)$; this correspondence is injective,
and we obtain an injection
${\mathcal{S}_\TOP(X)\too\RTOP(\sigma_TX)}$.
We call $\mathcal{T}_\TOP(X)=\RTOP(\sigma_TX)$ the
set of \emph{tangential invariants} of $X$. 
(Most texts consider normal invariants instead of tangential invariants.
Since we are working in the stable category, the difference is
merely the sign. Kirby-Siebenmann \cite{KS} use
tangential invariants.)
\end{Num}
Given a spherical fibration $\phi$ which admits a TOP reduction, it can be
shown that the abelian group $[X;\G/\TOP]$ acts regularly on
$\RTOP(\phi)$; thus there is a bijection of sets
\[
\mathcal{T}_\TOP(X)\cong[X;\G/\TOP].
\]
Now we can state the surgery classification of manifolds of a given
homotopy type. So let $X$ be a $1$-connected finite Poincar\'e duality
complex of formal dimension $n\geq 5$. There exists an abelian group
$P_n$ and a map
$\theta:\mathcal{T}_\TOP(X)\too P_n$, such that $\mathcal{S}_\TOP(X)$ is
precisely the preimage $\theta^{-1}(0)$; in other words, there is an exact
sequence of sets
\begin{diagram}[height=2em,width=6em]
\mathcal{S}_\TOP(X) & \rInto^{\text{inj}} &
\mathcal{T}_\TOP(X) & \too^\theta & P_n\\
&& \uTo & \ruDotsto \\
&& [X;\G/\TOP]
\end{diagram}
the \emph{surgery exact sequence}. The dotted arrow is in general
not a group homomorphism (and $\mathcal{T}_\TOP(X)$ has no canonical
group structure). If $n=4k\geq 8$ (the case we are interested in), then 
$P_n\cong\ZZ$ and $\theta$ is connected to the $\L$-genus and the
signature as
\[
\theta(\xi)=\frac18(\bra{\L_{4k}(\xi),[X]}-\mathrm{Sig}(X)).
\]
In other words, a bundle $\xi$ represents a homotopy manifold
structure for $X$
if and only if $\xi$ satisfies Hirzebruch's signature theorem.

Using these techniques,
the proof of Theorem \ref{ThmCompleteModels}
is accomplished by the following steps.

\paragraph{Step 1}
Every 1-connected Poincar\'e duality complex $P$ (of formal dimension
$2m\geq 5$) which has the same
homology as a projective plane is homotopy equivalent to one of
our models $M(\xi)$.

\paragraph{Step 2}
Let $\phi$ be a stable
spherical fibration over a finite, 1-connected CW-complex
$X$ with the property that $H_k(X)=0$ for $k\not\equiv0\pmod4$. Then
the Pontrjagin character $ph$ injects $\RTOP(\phi)$ into $H^\bullet(X;\QQ)$.
In particular, $ph:\mathcal{T}_\TOP(P)\too H^\bullet(P;\QQ)$ is
an injection.

\paragraph{Step 3}
For $m=4,8$,
we show that the elements of
$\mathcal{S}_\TOP(P)$ are completely determined by their  Pontrjagin classes
$p_m$, and that all possibilities for the $p_m$ are realized through our
models $M(\xi)$.

\paragraph{Step 4}
We determine $\mathrm{Aut}(P)$ and $\mathcal{M}_\TOP(P)$.

\begin{Num}
\label{HomoTopSetUp}
In the remainder of this section, we carry out steps 1 -- 4.
Let $P$ be a 1-connected Poincar\'e duality complex of formal dimension
$2m$ which is like a projective plane, for $m=4,8$.
So $H_k(P)\cong\ZZ$ for $m=0,1,2$.
We fix a map 
\[
s:\SS^m\too P
\]
representing a generator of
$\pi_m(P)\cong H_m(P)\cong\ZZ$.
Note also that $P$ has a preferred orientation $[P]$ -- the class dual
to $y_m^2$, where $y_m\in H^m(P)\cong\ZZ$ is any generator. Thus, any
homotopy equivalence automatically preserves fundamental classes.
By Wall \cite{WallCW} Proposition~4.1, we may assume that
\[
P=\SS^m\cup_\alpha e^{2m}
\]
is a $2$-cell complex, and that $m=4,8$, see \ref{CellStructure}.
\end{Num}
\paragraph{\em Proof for Step 1.}
This is just Theorem \ref{HomotopyTypesOfModels}. However, we
can say a bit more: the results in Section \ref{HomotopyOfModels}
show that the homotopy type of $P$ is uniquely determined by the
stable weak type of the spherical fibration $s^*\sigma_TP$, i.e.
by stable type of the pair of spherical fibrations
$\{s^*\sigma_TP,i^*s^*\sigma_TP\}$, where $i:\SS^m\rTo\SS^m$ is
any map of degree $-1$.

\paragraph{\em Proof for Step 2.}
Let $\phi$ be a stable spherical fibration over a connected
CW-complex $X$, with
stable classifying map $c:X\too\BG$. A (stable)
$\TOP$-reduction of $c$ is a lift $C$
\begin{diagram}[height=2em,width=4em]
&& \BTOP \\
&\ruTo^C &\dTo_{f_\G^\TOP} \\
X & \too^c & \BG,
\end{diagram}
with $f^\TOP_\G\circ C=c$. Two such lifts $C_0,C_1$ are equivalent if there
exists a homotopy $C:X\times[0,1]\too \BTOP$ with $c=f^\TOP_\G\circ C_t$
for all $t\in[0,1]$, i.e. if the homotopy is constant when projected to
$\BG$.
The set of equivalence classes of lifts of
$c$ is denoted $\RTOP(\phi)$.
If $c$ is the constant map (and thus $\phi=\underline 0$ is trivial),
we obtain a bijection
$\RTOP(\underline 0)\cong[X;\G/\TOP]$. For $\eta\in\RTOP(\phi)$ and
$\zeta\in\RTOP(\psi)$ we have $\eta\oplus\zeta\in\RTOP(\phi\oplus\psi)$.
This establishes a bijection
\begin{align*}
\RTOP(\underline 0) & \too \RTOP(\phi) \\
\eta & \mapstoo \eta\oplus\zeta
\end{align*}
see Wall \cite{Wallsurgery} Sec.~10, p.~113.
Thus, we can identify $\RTOP(\phi)$ with $[X;\G/\TOP]$, provided that
$\RTOP(\phi)\neq\emptyset$. Note however that $\RTOP(\phi)$ has
in general no natural group structure; rather, the abelian group
$[X;\G/\TOP]$ acts regularly on this set.

In general, different elements of $\RTOP(\phi)$ can be equivalent
when viewed as stable bundles. From the homotopy viewpoint, this
is due to the fact that two lifts $C,C'$ can be homotopic without
being fiber homotopic. We prove now that under certain conditions on $X$,
the map $\RTOP(\phi){\too}{[X,\BTOP]}$ is injective.
\begin{Prop}
\label{FiberHomotopyProp}
Let $X$ be a finite $1$-connected CW-complex. Assume that $H_k(X)=0$ for
all $k\not\equiv0\pmod 4$. Then the natural map
${[X;\G/\TOP]}{\too}{[X;\BTOP]}$
is injective.
\end{Prop}
Before we start with the proof, we note the following. By
Theorem \ref{Thetafinite}, we have isomorphisms
$\KOtilde(\SS^k)\otimes\QQ\rTo\KTOPtilde(\SS^k)\otimes\QQ$
for all $k\geq 0$. By a well-known comparison theorem
for half-exact cofunctors, see Dold \cite{Doldhalb} Ch.~7, or by
the Atiyah-Hirzebruch spectral sequence, see Hilton \cite{Hilton} Ch.~3,
this implies that
\[
\KOtilde(Y)\otimes\QQ\too^\cong\KTOPtilde(Y)\otimes\QQ
\]
is a natural isomorphism of homotopy functors for every finite 
connected CW-complex $Y$.
We combine this with the Pontrjagin character $ph$, which is 
(by the same comparison theorem for half-exact cofunctors and
by Bott Periodicity, see Section \ref{CharClass}) rationally
an isomorphism 
\[
\KOtilde(Y)\otimes\QQ\too^{ph}_\cong\widetilde H^{4\bullet}(Y;\QQ)
\]
to obtain an isomorphism
\[
\KTOPtilde(Y)\otimes\QQ\too^{ph}_\cong\widetilde H^{4\bullet}(Y;\QQ)
\]
which we also denote by $ph$.

\medskip
\emph{Proof of Proposition \ref{FiberHomotopyProp}}.
It clearly suffices to show that the natural map
\[
[X;\G/\TOP]\too\KTOPtilde(X)\otimes\QQ
\]
is an injection.
First we note that this is true in general for 
$X=\SS^{4k+t}$, for $t=-1,0,1$ (note that $\pi_{4k\pm1}(\G/\TOP)=0$,
cp.~Theorem \ref{SurgeryGroups}).
Both $[-;\G/\TOP]$ and $\KTOPtilde(-)\otimes\QQ$
are half-exact
cofunctors, so injectivity holds also for a wedge of spheres
$X=\bigvee_1^r\SS^{4k+t}$ (this is just the
additivity of half-exact cofunctors).

For a general complex $X$ as in the
claim of Proposition \ref{FiberHomotopyProp}, we proceed by induction.
By standard obstruction theory, we may assume that
$X$ is complex whose cells all have dimensions divisible by $4$,
i.e. that $X^{(4k)}=X^{(4k+1)}=X^{(4k+2)}=X^{(4k+3)}$ for all $k$,
see Wall \cite{WallCW} Proposition~4.1.
So suppose that $X=X^{(4k)}$, and that $A=X^{(4k-1)}=X^{(4k-4)}$.
The long exact sequence of the pair $(X,A)$ shows that $A$
also satisfies the hypothesis of Proposition \ref{FiberHomotopyProp},
and by induction, we may assume that the conclusion of the proposition 
holds for $A$.
Now consider the Puppe sequence
\[
\textstyle
\bigvee_r\SS^{4k-1}\too A \too X \too 
\bigvee_r\SS^{4k}\too SA \too\cdots
\]
Note that
$\KTOPtilde(SA)\otimes\QQ\cong\widetilde H^{4\bullet}(SA;\QQ)=0$, and
that similarly $\KTOPtilde(\bigvee_r\SS^{4k-1})\otimes\QQ=0$.
We thus obtain a diagram
{\scriptsize
\begin{diagram}
0 & \lTo &[A;\G/\TOP] & \lTo & [X;\G/\TOP] & \lTo &
{\textstyle{[\bigvee_r\SS^{4k};\G/\TOP]}} & \lTo
&[SA;\G/\TOP] \\
\dTo && \dTo^{\text{inj}} && \dTo && \dTo^{\text{inj}}
&& \dTo^{\text{inj}} \\
0& \lTo & \KTOPtilde(A)\otimes\QQ & \lTo &\KTOPtilde(X)\otimes\QQ
&\lTo &
{\textstyle\KTOPtilde(\bigvee_r\SS^{4k})\otimes\QQ}
&\lTo &
0 
\end{diagram}}%
and this implies by the Five-Lemma that
$[X;\G/\TOP]\too\KTOPtilde(X)\otimes\QQ$ is injective, see
eg. Eilenberg-Steenrod \cite{ES} Lemma 4.4.
This finishes the proof of Proposition \ref{FiberHomotopyProp}.\qed

\begin{Cor}
\label{CorStructureSet}
Let $\phi$ be a spherical fibration over a finite 1-connected
CW-complex $X$, with $H_k(X)=0$ for all $k\not\equiv 0\pmod4$,
and assume that
$\RTOP(\phi)\neq\emptyset$. Then $\RTOP(\phi)$ injects into
$\KTOPtilde(X)\otimes\QQ$ and, via the Pontrjagin character, into
$\widetilde H^{4\bullet}(X;\QQ)$.

\proof
Let $\zeta\in \RTOP(\phi)$ be a stable bundle. By Proposition
\ref{FiberHomotopyProp}, the Pontrjagin character injects
$\RTOP(\underline 0)$ into $\widetilde H^{4\bullet}(X;\QQ)$.
The diagram {\small
\begin{diagram}[width=5em]
\RTOP(\underline 0) & \too^{[\eta\mapstoo\eta+\zeta]}_{\text{bij}} &
\RTOP(\phi)\\
\dTo_{ph}^{\text{inj}} && \dTo_{ph} \\
\widetilde H^{4\bullet}(X;\QQ) & 
\too^{[x\mapstoo x+ph(\zeta)]}_{\text{bij}} &
\widetilde H^{4\bullet}(X;\QQ)
\end{diagram}}%
commutes, and the claim follows.
\qed
\end{Cor}
\begin{Cor}
Let $X$ be a finite 
1-connected Poincar\'e duality complex, and assume that
$H_k(X)=0$ for all $k\not\equiv0\pmod 4$. Then the Pontrjagin character
$ph$ injects the set $\mathcal{T}_\TOP(X)$of tangential invariants
into $\widetilde H^{4\bullet}(X;\QQ)$.
\qed
\end{Cor}
This corollary applies in particular to our 
Poincar\'e duality complex $P$ of formal dimension $2m$, for $m=4,8$.
Note also the following. If $\eta,\zeta\in \RTOP(\phi)$ are elements
with the same total Pontrjagin class, $p(\eta)=p(\zeta)$, then
clearly $ph(\eta)=ph(\zeta)$. Therefore, the total Pontrjagin class
$p$ induces also an injection of $\RTOP(\phi)$ into
$\widetilde H^{4\bullet}(X;\QQ)$.

\paragraph{\em Proof for Step 3.}
We use the same symbols $P$, $y_m,y_{2m}$, $s$ as in \ref{HomoTopSetUp}.
Since $\mathrm{Sig}(P)=1$, a stable TOP-bundle reduction
$\zeta\in\mathcal{T}_\TOP(P)$
represents a homotopy manifold structure for $P$ if and only if
\[
\L_{2m}(\zeta)=y_{2m}^2.
\]
Also, we see from the formula for the $\L$-genus \ref{SpecificValues}
that
there are rational numbers $c_m,d_m$ (depending only on $m$)
such that
\[
p_{2m}(\eta)=c_m\L_{2m}(\eta)+d_m p_m(\eta)^2
\]
for any stable bundle $\eta$ over $P$.
\begin{Lem}
\label{LemmaTheRestrictionClassifies}
The map $\mathcal{T}_\TOP(P)\too H^m(P;\QQ)$,
$\zeta\mapstoo p_m(\zeta)$ is an injection when restricted to
$\mathcal{S}_\TOP(P)\subseteq\mathcal{T}_\TOP(P)$.

\proof
Let $\eta,\zeta$ be elements in $\mathcal{T}_\TOP(X)$ representing
homotopy manifold structures, so $\L_{2m}(\eta)=y_m^2=\L_{2m}(\zeta)$.
If $p_m(\eta)=p_m(\zeta)$, then $p_{2m}(\eta)=p_{2m}(\zeta)$
by the formula above.
Thus $p(\eta)=p(\zeta)$, whence $ph(\eta)=ph(\zeta)$, and
therefore $\eta=\zeta$ by Corollary \ref{CorStructureSet}.
\qed
\end{Lem}
Now we prove that our set of model manifolds realizes all elements
in $\mathcal{S}_\TOP(P)$. Let $\phi=s^*\sigma_TP$.
\begin{Lem}
\label{CanFindEuler1}
Let $\zeta\in\RTOP(\phi)$. Then there exists an
$\RR^m$-bundle $\xi$ over $\SS^m$ with absolute Euler number $|e|=1$
which is stably equivalent to $\zeta$.

\proof
We have $w_m(\zeta)=x$ mod $2$ (because $w_m(M_{r,s})=y_m$ mod $2$ for any
model). From the split exact sequence in Lemma \ref{SplittingLemma}
we see that we can find an $\RR^m$-bundle $\xi$ over $\SS^m$ with
$w_m(\xi)=x$ mod $2$, and hence with any odd absolute Euler number.
\qed
\end{Lem}
\begin{Cor}
For $m=4,8$, all elements of $\mathcal{S}_\TOP(X)$ are realized as
model manifolds $M(\xi)$.

\proof
Given a stable bundle $\eta\in\mathcal{T}_\TOP(P)$ representing
a homotopy manifold structure $M\too X$ in $\mathcal{S}_\TOP(X)$,
we can find by Lemma \ref{CanFindEuler1} an oriented $\RR^m$-bundle
$\xi$ over $\SS^m$ which is stably equivalent to $s^*\eta$,
with absolute Euler number $|e|=1$.
Then the model manifold $M(\xi)$ is homotopy equivalent to $P$ by
the remark in Step 1, because the homotopy types of $P$ and
$M(\xi)$ are
determined by $s^*\sigma_TP$ and $s_0^*\sigma_TM(\xi)$,
respectively.
Composing the homotopy equivalence $M(\xi)\simeq P$
with a self-homeomorphism of $M(\xi)$ induced by a homeomorphism
of degree $-1$ of $\SS^m$, if necessary, we obtain a homotopy
commutative diagram
\begin{diagram}[size=2em]
\SS^m &&\too^s&& P\\
&\rdTo_{s_0} && \ldTo^f_\simeq \\
&& M(\xi).
\end{diagram}
By Lemma \ref{LemmaTheRestrictionClassifies}, $M(\xi)$ is precisely
the homotopy manifold structure on $P$ represented by $\eta$, because
$f^\bullet p_m(M)=p_m(\eta)$, so $M\cong M(\xi)$.
\qed
\end{Cor}
This finishes Step 3.

\paragraph{\em Proof for Step 4.}
We proved already in Proposition \ref{HomeoProp}
that $M(\xi)\cong M(\xi')$ if and only if
$\xi$ and $\xi'$ are weakly equivalent.
The group
$\mathrm{Aut}(M(\xi))$ is cyclic of order two by Lemma \ref{AutX},
and this finishes the classification for $m=4,8$.

\section{Appendix: homotopy classification}

\emph{In this section, all maps and homotopies
are assumed to preserve base points.}

\begin{Num}
\label{CellStructure}
Suppose that $X$ is a $1$-connected Poincar\'e duality complex 
of formal dimension $n$, with $H_\bullet(X)\cong\ZZ^3$
(cp. \ref{PDCDef}).
Let $\mu\in H_n(X)$ denote the fundamental class. We have
$H_0(X)\cong\ZZ\cong H_n(X)$, so $H_m(X)\cong\ZZ$ for some number
$1<m<n$. By the universal coefficient theorem, see Spanier
\cite{Span} Ch.~5.5 Theorem 4, we have $H^j(X)\cong H_j(X)$ for all $j$.
From Poincar\'e duality, we see that $n=2m$, and that the map
\[
H^m(X)\otimes H^m(X)\too\ZZ,\qquad
u\otimes v\mapstoo\bra{u{\smile}v,\mu}
\]
is a duality pairing. Thus $m$ is even, and
$H^\bullet(X)\cong\ZZ[y_m]/(y_m^3)$, for some generator $y_m\in H^m(X)$.
By Wall \cite{WallCW} Proposition~4.1, the CW-complex $X$ is homotopy
equivalent to a $2$-cell complex, 
\[
X\simeq X_\alpha=\SS^m\cup_\alpha e^{2m},
\]
for some attaching map $\alpha:\SS^{2m-1}\rTo\SS^m$.
By Adams-Atiyah \cite{AA} Theorem~A, this implies that $m=2,4,8$.
Note also that $X_\alpha$ has a preferred orientation $\mu$, the
dual of $y_m^2$.
\end{Num}
\begin{Num}
We wish to determine the number of homotopy types of such complexes $X$.
The structure of the cohomology ring of $X_\alpha$ implies that
$\alpha$, viewed as an element of $\pi_{2m-1}(\SS^m)$, has
Hopf invariant $h(\alpha)=\pm1$,
see Adams-Atiyah \cite{AA} and Husemoller \cite{Hus}
Ch.~20.10. Note also that the homotopy type of $X_\alpha$ is not
changed if $\alpha$ is replaced by a map homotopic to $\alpha$,
see Milnor \cite{MilnorMorse} Lemma 3.6.
Let $\mathbf W_m$ denote the set of all homotopy types of
$1$-connected CW-complexes as above. Let
$\mathcal H_m^{\pm1}\subseteq\pi_{2m-1}(\SS^m)$ denote the set of
all elements of Hopf invariant $\pm 1$. We thus have a well-defined 
surjection
\[
\pi_{2m-1}(\SS^m)\supseteq \mathcal H_m^{\pm1}\rTo\mathbf W_m,
\]
sending $\alpha\in\mathcal H_m^{\pm1}$ to the homotopy type of $X_\alpha$.
By Adams-Atiyah \cite{AA} Theorem~A, the set $\mathbf W_m$ is
empty unless $m=2,4,8$.

Next we note the following. If $c:\SS^{2m-1}\rTo\SS^{2m-1}$ is an
involution of degree $-1$, then there is a homeomorphism
$X_\alpha\cong X_{\alpha \circ c}$. Since $\alpha\circ c$ represents
$-\alpha\in \pi_{2m-1}(\SS^m)$, we have a homotopy equivalence
\[
X_\alpha\simeq X_{-\alpha}.
\]
Each element in $\mathbf W_m$ is thus represented by a element
$\alpha\in\pi_{2m-1}(\SS^m)$ with Hopf invariant $h(\alpha)=1$,
i.e. $\mathcal H_m^+$ surjects onto $\mathbf W_m$.
\end{Num}
For $m=2$, we are done: $\pi_3(\SS^2)\cong\ZZ$, see Toda \cite{Toda}
p.~186, so $\mathcal H_2^{\pm}$ has exactly two elements, and
$\mathbf W_2$ consists of precisely one homotopy type, the complex
projective plane $\CP^2$.
\begin{Lem}
There is exactly one homotopy type in $\mathbf W_2$.
\qed
\end{Lem}
\begin{Num}
It remains to consider the cases $m=4,8$.
Similarly as above, if $g:\SS^m\rTo\SS^m$ is an involution of degree $-1$,
then there is a homeomorphism
$X_\alpha\cong X_{g\circ\alpha}$, so
\[
X_\alpha\simeq X_{g_\#(\alpha)}.
\]
The homotopy equivalences 
\[
X_\alpha\simeq X_{-\alpha}\simeq X_{g_\#(\alpha)}\simeq X_{-g_\#(\alpha)}
\]
are in fact the only homotopy equivalences which occur between these
$2$-cell complexes.
For suppose that $f:X_\alpha\rTo^\simeq X_\beta$ is a homotopy equivalence.
We may assume that $f$ is a cellular map, see Whitehead \cite{Whi}
Ch.~II Theorem 4.5, 
\[
f:(X_\alpha,\SS^m)\rTo(X_\beta,\SS^m), 
\]
and then
$f$ restricts to a map of degree $\pm 1$ on the $m$-skeleton $\SS^m$.
Also, we have the Hurewicz isomorphism
\[
H_{2m}(\SS^{2m})=H_{2m}(X_\alpha/\SS^m)\cong
H_{2m}(X_\alpha,\SS^m)\cong \pi_{2m}(X_\alpha,\SS^m)\cong \ZZ;
\]
a canonical generator $\hat\alpha$
of this group $\pi_{2m}(X_\alpha,\SS^m)$ is given by the attaching map,
\[
(e^{2m},\SS^{2m-1})\rTo^{\hat\alpha}(X_\alpha,\SS^m)
\quad\text{ with }\quad
\partial\hat\alpha=\alpha.
\]
From the homotopy exact sequence
\begin{diagram}[height=2em]
{}&\rTo & \pi_{2m}(\SS^m) & \rTo & \pi_{2m}(X_\alpha) & \rTo
&\pi_{2m}(X_\alpha,\SS^m) & \rTo^\partial & \pi_{2m-1}(\SS^m) & \rTo&{}\\
&&\dTo && \dTo&&\dTo&&\dTo\\
{}&\rTo & \pi_{2m}(\SS^m) & \rTo & \pi_{2m}(X_\beta) & \rTo
&\pi_{2m}(X_\beta,\SS^m) & \rTo^\partial & \pi_{2m-1}(\SS^m) & \rTo&{}
\end{diagram}
and the  Five-Lemma,
we see that $f_\#(\hat\alpha)=\pm\hat\beta$, and so $f_\#(\alpha)=\pm\beta$.
We have seen above that we may replace $X_\beta$ by $X_{-\beta}$;
thus, we may assume that $f_\#(\alpha)=\beta$. If $f$ restricts to
a map of degree $1$ on $\SS^m$, then $f_\#(\alpha)=\alpha=\beta$.
\end{Num}
\begin{Lem}
\label{AutX}
The group of self-equivalences $\mathrm{Aut}(X_\alpha)$
is cyclic of order $2$;
it coincides with the group of graded ring automorphisms
of the cohomology ring $\ZZ[y_m]/(y_m^3)$.
\qed
\end{Lem}
\begin{Num}
So the remaining problem is to determine the relation between
$\alpha$ and $g_\#(\alpha)$, where $g:\SS^m\rTo\SS^m$ is a map
of degree $-1$.
Towards this end, we consider the $EHP$-sequence of
$\SS^n$ for the values $n=m-1,m$,
\[
\rTo\pi_k(\SS^n)\rTo^E
\pi_{k+1}(\SS^{n+1})\rTo^H\pi_{k+1}(\SS^{2n+1})\rTo^P\pi_{k-1}(\SS^n)\rTo^E
\]
This sequence is exact for $k\leq3n-2$,
see Whitehead \cite{Whi} Ch.~XII Theorem 2.2.
Here, $E$ is the suspension and $H$ is the
generalized Hopf invariant. Let $\iota_j=\id_{\SS^j}$ denote the canonical
generator of $\pi_j(\SS^j)$. For $\rho\in\pi_{2n+1}(\SS^n)$, one
has $H(\rho)=h(\rho)\cdot\iota_{2n+1}$; see Whitehead
\cite{Whi} for a comparison
between the various definitions of Hopf invariants.
\end{Num}
\begin{Num}
Consider the diagram
{\small
\begin{diagram}
&&&& 0 \\
&&&& \dTo \\
&&\pi_{2m-1}(\SS^{2m-1})&\too^{E^2}_\isom& \pi_{2m+1}(\SS^{2m+1})\\
&&&& \dTo_P & \rdTo^{\cdot 2} \\
0 & \too & 
\pi_{2m-2}(\SS^{m-1}) & \too^E &
\pi_{2m-1}(\SS^{m}) & \too^H &
\pi_{2m-1}(\SS^{2m-1}) & \too & 0 \\
&&& \rdTo & \dTo_E \\
&&&& \pi_{2m}(\SS^{m+1}) \\
&&&& \dTo \\
&&&& 0\rlap{.}
\end{diagram}}%
The middle row is split and
short exact: $\pi_{2m-1}(\SS^{2m-1})$ is
infinite cyclic (whence the splitting)
and the Hopf invariant $H$ is by assumption
onto. The map $P:\pi_{2m}(\SS^{2m-1})\rTo\pi_{2m-2}(\SS^{m-1})$ can be
characterized by $PE^2(\rho)=[\iota_{m-1},\iota_{m-1}]\circ\rho$,
see Whitehead \cite{Whi} Ch.~XII Theorem 2.4. But $\SS^{m-1}$ is an
$H$-space for $m=4,8$, so $[\iota_{m-1},\iota_{m-1}]=0$, see
Whitehead \cite{Whi} Ch.~X Corollary 7.8. From the
$EHP$ sequence, we see that
$\pi_{2m-2}(\SS^{m-1})\too^E\pi_{2m-1}(\SS^{m})$ is an injection.

The middle column is also short exact: from Freudenthal's
Suspension Theorem we have that
$E:\pi_{2m-1}(\SS^m)\rTo\pi_{2m}(\SS^{m+1})$
is an epimorphism, see Whitehead \cite{Whi} Ch.~VII Theorem 7.13.
To see that $P:\pi_{2m+1}(\SS^{2m+1})\rTo\pi_{2m-1}(\SS^m)$ is
injective, note that $PE^2(\iota_{2m-1})=[\iota_m,\iota_m]$ by
Whitehead \cite{Whi} Ch.~XII Theorem 2.4. But
\[
H([\iota_m,\iota_m])=2\iota_{2m-1},
\]
see Whitehead \cite{Whi} Ch.~XI Theorem 2.5. Thus, $P$ is injective
on this infinite cyclic group.
\end{Num}
\begin{Num}
So suppose that $g:\SS^m\rTo\SS^m$ has degree $-1$.
Then 
\[
g_\#([\iota_m,\iota_m])=
[g_\#(\iota_m),g_\#(\iota_m)]=[-\iota_m,-\iota_m]=[\iota_m,\iota_m],
\]
whence 
\[
H(g_\#(\xi))=H(\xi)\text{ for all }\xi\in\pi_{2m-1}(\SS^m).
\]
Now let $\rho\in\pi_{2m-2}(\SS^{m-1})$. This group is finite by Serre's
finiteness result for odd spheres,
see Spanier \cite{Span} Ch.~9.7 Theorem 7; via $E$, we can identify it
with the torsion group of $\pi_{2m-1}(\SS^m)$.
The double suspension $E^2$ injects this group into the stable
group $\pi_{2m}(\SS^{m+1})=\pi_{m-1}^s(\SS^0)$. In the graded ring
$\pi_\bullet^s(\SS^0)$, composition is commutative, see Whitehead
\cite{Whi} Ch.~XII; thus we have $E(g_\#(E\rho))=-E^2(\rho)$,
whence $g_\#(E(\rho))=-E(\rho)$. The involution $g_\#$ thus changes
the signs of the elements of the torsion group of $\pi_{2m-1}(\SS^m)$.

The stable groups $\pi_{m-1}^s(\SS^0)$ are cyclic of order $24$ and $240$,
for $m=4,8$, see Toda \cite{Toda} p.~186. The image of the
double suspension $E^2$ is also
cyclic and has index $2$ in this group.
Thus, if
$\xi\in\pi_{2m-1}(\SS^m)$ is an element with Hopf invariant $1$,
then the element $E(\xi)$ together with the cyclic group
$E^2(\pi_{2m-2}(\SS^{m-1}))$ generates
$\pi_{2m+2}(\SS^{m+1})\cong\pi_{m-1}^s(\SS^0)$.
Therefore, we can find an element
$\eta_m\in\pi_{2m-1}(\SS^m)$ with Hopf invariant $h(\eta_m)=1$
whose suspension $E(\eta_m)$ generates
$\pi_{2m}(\SS^{m+1})$.
Let $\delta_m=2\eta_m-[\iota_m,\iota_m]$; then
$E(\delta_m)=2\cdot E(\eta_m)$
generates the cyclic group $E^2(\pi_{2m-2}(\SS^{m-1}))$
(the suspension of the Whitehead product $[\iota_m,\iota_m]$
vanishes, see Whitehead \cite{Whi} Ch.~X Theorem 8.20), so
$\delta_m$ generates the torsion group of $\pi_{2m-1}(\SS^m)$.
Put $g_\#(\eta_m)=\eta_m+r\cdot\delta_m$. Now
\[
-\delta_m=g_\#(\delta_m)=g_\#(2\eta_m-[\iota_m,\iota_m])=
2\eta_m+2r\cdot\delta_m-[\iota_m,\iota_m]=(1+2r)\cdot\delta_m,
\]
whence $2(1+r)\delta_m=0$. This leaves two possibilities for $r$,
for $m=4,8$. But due to the commutativity of
$\pi_\bullet^s(\SS^0)$, we have also
$E(g_\#(\eta_m))=-E(\eta_m)=E(\eta_m)+r\cdot E(\delta_m)=(1+2r)E(\eta_m)$,
which implies that $r=-1$, i.e.~that 
\[
g_\#(\eta_m)=\eta_m-\delta_m.
\]
We have proved the following result.
\end{Num}
\begin{Prop}
The number of homotopy types in $\mathbf W_m$ is $1,6,60$, for
$m=2,4,8$.
\qed
\end{Prop}
This is a complete homotopy classification of manifolds and
complexes which are
like projective planes, and also the end of this paper.

\bibliographystyle{plain}

\begin{thebibliography}{WWW}

\bibitem{AdamsJ-IV}
J.F. Adams,
On the groups $J(X)$--IV,
Topology 5 (1966) 21--71.

\bibitem{AA}
J.F. Adams and M.F. Atiyah,
$K$-theory and the Hopf invariant,
Quart J. Math. Oxford 17 (1966) 31--38.

\bibitem{AH}
M.F. Atiyah and F. Hirzebruch,
Vector bundles and homogeneous spaces,
Proc. Symp. Pure Math. 3 (1961) 7--38.

\bibitem{Bott1}
R. Bott,
The stable homotopy of the classical groups,
Ann. Math. 70 (1959) 313--227.

\bibitem{Bott2}
R. Bott,
The periodicity theorem for the classical groups and some
of its applications,
Adv. Math. 4 (1970) 353--411.

\bibitem{BK}
U. Brehm and W. K\"uhnel,
$15$-vertex triangulations of an $8$-manifold,
Math. Ann. 294 (1992) 167--193.

\bibitem{Bs} 
S. Breitsprecher, 
Projektive Ebenen, die Mannigfaltigkeiten sind, 
Math. Z. 121 (1971) 157--174.

\bibitem{Brumfiel}
G. Brumfiel,
On the homotopy groups of $\BPL$ and $\PL/\O$,
Ann. Math. 88 (1968) 291--311.

\bibitem{Buch}
T. Buchanan,
Zur Topologie der projektiven Ebenen \"uber reellen
Divisionsalgebren,
Geom. Ded. 8 (1979) 383--393.

\bibitem{CE} 
H. Cartan and S. Eilenberg,
{\em Homological algebra},
Princeton University Press, Princeton, New Jersey (1956).

\bibitem{Dold}
A. Dold,
Partitions of unity in the theory of fibrations,
Ann. Math. 78 (1963) 223-255.

\bibitem{Doldhalb}
A. Dold,
\emph{Halbexakte Homotopiefunktoren},
Springer LNM 12,
Springer-Verlag (1966).

\bibitem{DoldAT}
A. Dold,
{\em Lectures on algebraic topology,}
Springer Verlag, Berlin Heidelberg New York (1972).

\bibitem{EK} 
J. Eells and N.H. Kuiper, 
Manifolds which are like projective planes, 
Publ. Math. I.H.E.S. 14 (1962) 5--46.

\bibitem{EK2} 
J. Eells and N.H. Kuiper, 
An invariant for certain smooth manifolds,
Ann. Mat. Pura Appl. 60 (1962) 93--110.

\bibitem{ES} 
S. Eilenberg and N. Steenrod,
{\em Foundations of algebraic topology},
Princeton University Press, Princeton, New Jersey (1952).

\bibitem{EngDim}
R. Engelking,
{\em Theory of dimensions, finite and infinite,}
Heldermann, Lemgo (1995).

\bibitem{FFG} 
A.T. Fomenko, D.B. Fuchs, and V.L. Gutenmacher,
{\em Homotopic topology,}
Akad\'emiai Kiad\'o, Budapest (1986)

\bibitem{Fr} 
M.H. Freedman,
The topology of four-dimensional manifolds,
J. Differential Geometry 17 (1982) 357--453.

\bibitem{FQ}
M.H. Freedman and F. Quinn,
\emph{Topology of $4$-manifolds},
Princeton Univ. Press, Princeton (1990).

\bibitem{Han}
O. Hanner,
Some theorems on absolute neighborhood retracts,
Arkiv Mat.~1 (1950) 389--408.

\bibitem{Hilton}
P. Hilton,
\emph{General cohomology theory and $K$-theory},
Cambridge Univ. Press, Cambridge (1971)

\bibitem{Hir}
M.W. Hirsch,
Obstruction theories for smoothing manifolds and maps,
Bull. Amer. Math. Soc. 69 (1963) 352--356.

\bibitem{HirNeue}
F. Hirzebruch,
\emph{Neue topologische Methoden in der algebraischen Geometrie,}
zweite erg. Auflage,
Springer-Verlag Berlin G\"ottingen Heidelberg (1962).

\bibitem{Hirz}
F. Hirzebruch,
A Riemann-Roch theorem for differentiable manifolds,
S\'em. Bourbaki Exp. 177, Paris (1959);
reprinted in: {\em F. Hirzebruch, Ges. Abh. Band I,}
pp. 481--495,
Springer Verlag, Berlin Heidelberg New York (1987)

\bibitem{Ho}
P. Holm,
The microbundle representation theorem,
Acta math. 117 (1967) 191--213.

\bibitem{Hu} 
S.-T. Hu,
{\em Homotopy theory,}
Academic Press (1959).

\bibitem{HuRetr} 
S.-T. Hu,
{\em Theory of retracts,}
Wayne State University Press (1965).

\bibitem{HuWa}
W. Hurewicz and H. Wallman,
{\em Dimension theory,}
Princeton University Press, Princeton (1948).

\bibitem{Hus} 
D. Husemoller,
{\em Fibre bundles, 3rd ed.,}
Springer Verlag (1994)

\bibitem{Kahn}
P. Kahn,
A note on topological Pontrjagin classes and the Hirzebruch
index formula,
Illinois J. Math. 16 (1972) 234--256.

\bibitem{KM}
M. Kervaire and J. Milnor,
Groups of homotopy spheres, I,
Ann. Math. 77 (1963) 504--537

\bibitem{KS} 
R.C. Kirby and L.C. Siebenmann, 
{\em Foundational essays on topological manifolds, smoothings, 
and triangulations,}
Princeton University Press, 
Princeton, New Jersey (1977).

\bibitem{Klein}
J.R. Klein,
Poincar\'e duality spaces, pp.~135--166, in:
\emph{Surveys on surgery theory, Vol. I},
S. Cappell, A. Ranicki and J. Rosenberg eds.,
Princeton Univ. Press (2000).

\bibitem{Kn} 
N. Knarr, 
The nonexistence of certain topological polygons,
Forum Math. 2 (1990) 603--612.

\bibitem{KrSmooth} 
L. Kramer, 
The topology of smooth projective planes,
Arch. Math. 63 (1994) 85--91.

\bibitem{KrP=L} 
L. Kramer, 
The point and line space of a compact projective plane are
homeomorphic,
Math. Z. 234 (2000) 83--102.

\bibitem{Kreck}
M. Kreck,
A guide to the classification of manifolds, pp.~121--134, in:
\emph{Surveys on surgery theory, Vol. I},
S. Cappell, A. Ranicki and J. Rosenberg eds.,
Princeton Univ. Press (2000).

\bibitem{LM} 
H.B. Lawson and M.-L. Michelsohn,
{\em Spin geometry},
Princeton University Press, Princeton, New Jersey (1989).

\bibitem{Lo} 
R. L\"owen, 
Topology and dimension of stable planes: On a conjecture of H. 
Freudenthal, 
J. reine angew. Math. 343 (1983) 108--122.

\bibitem{MM} 
I. Madsen and R.J. Milgram, 
{\em The classifying spaces for surgery and cobordism of manifolds,}
Princeton University Press, Princeton, New Jersey (1979).

\bibitem{MC} 
J. McCleary,
{\em A user's guide to spectral sequences, 2nd ed.},
Cambridge University Press, Cambridge (2001).

\bibitem{MiMicro}
J. Milnor,
Microbundles--I.
Topology 3 Suppl.~1 (1964) 53--80.

\bibitem{MilSpher} 
J.W. Milnor, 
On characteristic classes for spherical fibre spaces, 
Comm. Math. Helv. 43 (1968) 51--77.

\bibitem{MilnorMorse} 
J.W. Milnor and J.D. Stasheff, 
{\em Morse Theory,}
Princeton University Press, Princeton, New Jersey (1963).

\bibitem{MS} 
J.W. Milnor and J.D. Stasheff, 
{\em Characteristic classes,}
Princeton University Press, Princeton, New Jersey (1974).

\bibitem{New}
H.M.A. Newman,
The engulfing theorem for topological manifolds,
Ann. Math. 84 (1966) 555--571.

\bibitem{Okabe}
T. Okabe,
The existence of topological Morse function and its
application to topological $h$-cobordism theorem,
J. Fac. Sci. Univ. Tokyo Sect. IA Math 18 (1971) 23--35.

\bibitem{Qui}
F. Quinn,
Ends of maps, I,
Ann. Math. 110 (1979) 275--331.

\bibitem{Rud}
Y.B. Rudyak,
{\em On Thom spectra, orientability, and cobordism},
Springer Verlag (1998)

\bibitem{Ros}
J. Rosenberg,
$C^*$-algebras, positive scalar curvature, and the Novikov
conjecture -- III,
Topology 25 (1986) 319--336.

\bibitem{CPP} 
H. Salzmann, D. Betten, T. Grundh\"ofer, H. H\"ahl, 
R. L\"owen, and M. Stroppel, 
{\em Compact projective planes,}
de Gruyter Expositions in Mathematics 21, Berlin (1995).

\bibitem{Shim}
N. Shimada,
Differentiable structures in the $15$-sphere and Pontrjagin
classes of certain manifolds,
Nagoya Math. J. 12 (1957) 59--69.

\bibitem{Sma}
S. Smale,
Generalized Poincar\'e's conjecture in dimensions greater than
four,
Ann. Math. 74 (1961) 391--406.

\bibitem{Span} 
E.H. Spanier,
{\em Algebraic topology,}
Springer (1966).

\bibitem{Spivak}
M. Spivak,
Spaces satisfying Poincar\'e duality,
Topology 6 (1967) 77--101.

\bibitem{Steenrod}
N.E. Steenrod,
\emph{Topology of fibre bundles},
Princeton Univ. Press, Princeton (1951).

\bibitem{Stolz}
S. Stolz,
\emph{Hochzusammenh\"angende Mannigfaltigkeiten und ihre R\"ander},
Springer LNM 1116, Springer Verlag (1985).

\bibitem{Stolzpscm}
S. Stolz,
Simply connected manifolds of positive scalar curvature,
Ann. Math. 136 (1992) 511--540.
\bibitem{Toda}
H. Toda,
\emph{Composition methods in homotopy groups of spheres},
Princeton Univ. Press, Princeton (1962).

\bibitem{Var}
K. Varadarajan,
Span and stably trivial bundles,
Pacific J. Math. 60 (1975) 277--287.

\bibitem{Walln-1}
C.T.C. Wall,
Classification of $(n-1)$-connected $2n$-manifolds,
Ann. Math. 75 (1962) 163--189.

\bibitem{WallCW} 
C.T.C. Wall,
Finiteness conditions for CW-complexes,
Ann. Math. 81 (1965) 56--69.

\bibitem{Wallsurgery} 
C.T.C. Wall,
\emph{Surgery on compact manifolds, 2nd ed.},
Mathematical surveys and monographs vol. 69,
Amer. Math. Soc. (1999).

\bibitem{Web}
C. Weber,
Quelques th\'eoremes bien connus sur les A.N.R et les
C.W. complexes.
Enseign. Math. 13 (1968) 211--222.

\bibitem{Whi}
G.W. Whitehead,
\emph{Elements of homotopy theory,}
Springer-Verlag New York Berlin Heidelberg (1978).

\bibitem{Wil}
R.E. Williamson, Jr.,
Cobordism of combinatorial manifolds,
Ann. Math. 83 (1966) 1--33.

\end{thebibliography}

\bigskip
\raggedright
Linus Kramer \\
Mathematisches Institut \\
Universit\"at W\"urzburg \\
Am Hubland \\
D-97074 W\"urzburg \\
Germany \\
\makeatletter
{\tt kramer{@}mathematik.uni-wuerzburg.de}

\end{document}